\documentclass[12pt]{amsart}
\usepackage{graphicx}
\newtheorem{theorem}{Theorem}[section]
\newtheorem{corollary}[theorem]{Corollary}
\newtheorem{remark}[theorem]{Remark}
\newtheorem{lemma}[theorem]{Lemma}
\newtheorem{proposition}[theorem]{Proposition}

\numberwithin{equation}{section}
\author{Xiaoli Han, Jiayu Li, Jun Sun}

\address{Xiaoli Han, Department of Mathematical Sciences, Tsinghua University \\ Beijing 100084, P. R. of China.}
\email{xlhan@math.tsinghua.edu.cn}

\address{Jiayu Li, School of Mathematical Sciences, University of Science and Technology of China Hefei 230026 \\ AMSS CAS Beijing 100190, P. R. China}
\email{jiayuli@ustc.edu.cn}

\address{Jun Sun, School of Mathematical Sciences\\
Wuhan University\\ Wuhan 430072, P. R. of China.}
\email{sunjun@whu.edu.cn}

\thanks {The research was supported by the National Natural Science Foundation of China,  No.11131007, No.11471014., No. 11401440.}

\begin{document}

\title[The deformation of symplectic critical surfaces]
{The deformation of symplectic critical surfaces in a K\"ahler
surface-I}

\begin{abstract}
In this paper we derive the Euler-Lagrange equation of the
functional $L_\beta=\int_\Sigma\frac{1}{\cos^\beta\alpha}d\mu,
~~\beta\neq -1$ in the class of symplectic surfaces. It is
$\cos^3\alpha {\bf{H}}=\beta(J(J\nabla\cos\alpha)^\top)^\bot$,
which is an elliptic equation when $\beta\geq 0$. We call such a
surface a $\beta$-symplectic critical surface. We first study the
properties for each fixed $\beta$-symplectic critical surface and
then prove that the set of $\beta$ where there is a stable
$\beta$-symplectic critical surface is open. We believe it should
be also closed. As a precise example, we study rotationally
symmetric $\beta$-symplectic critical surfaces in ${\mathbb C}^2$
carefully .
\end{abstract}
\maketitle

\section{Introduction}

\allowdisplaybreaks

Suppose that $M$ is a K\"ahler surface. Let $\omega$ be the
K\"ahler form on $M$ and let $J$ be a complex structure compatible
with $\omega$. The Riemannian metric $\langle,\rangle$ on $M$ is
defined by
$$
\langle U,V \rangle =\omega(U,JV).
$$
For a compact oriented real surface $\Sigma$ which is smoothly
immersed in $M$, one defines, following \cite {CW}, the K\"ahler
angle $\alpha$ of $\Sigma$ in $M$ by
\begin{equation}\label{e1}\omega|_\Sigma=\cos\alpha d\mu_\Sigma\end{equation} where $d\mu_\Sigma$
is the area element of $\Sigma$ of the induced metric from
$\langle,\rangle$. We say that $\Sigma$ is a holomorphic curve if
$\cos\alpha \equiv 1$, $\Sigma$ is a Lagrangian surface if
$\cos\alpha \equiv 0$ and $\Sigma$ is a symplectic surface if
$\cos\alpha > 0$.

In \cite{HL} we consider the functional
$$L=\int_\Sigma\frac{1}{\cos\alpha}d\mu.$$
The Euler-Lagrange equation of this functional is
$$\cos^3\alpha {\bf{H}}=(J(J\nabla\cos\alpha)^\top)^\bot.$$ We call
such a surface a symplectic critical surface. We study the
properties of the symplectic critical surfaces. In this paper, we
consider a sequence of functionals
$$L_\beta=\int_\Sigma\frac{1}{\cos^\beta\alpha}d\mu.$$

The critical point of the functionals $L_\beta$ in the class of
symplectic surfaces in a K\"ahler surface is called a
{\em $\beta$-symplectic critical surface}. We first calculate the
Euler-Lagrange equation of $L_\beta$.

\begin{theorem}\label{crie} Let $M$ be a K\"ahler surface.
The Euler-Lagrange equation of the functional $L_\beta$
($\beta\neq -1$) is \begin{equation}\label{me}\cos^3\alpha
{\bf{H}}-\beta(J(J\nabla\cos\alpha)^\top)^\bot=0, \end{equation}
where $H$ is the mean curvature vector of $\Sigma$ in $M$, and
$()^\top$ means tangential components of $()$, $()^\bot$ means the
normal components of $()$.
\end{theorem}

Similar to the case that $\beta=1$ in \cite{HL}, we can check that it
is an elliptic equation module tangential diffeomorphisms if $\beta\geq 0$.

The equation (\ref{me}) seems quit interesting. It is the minimal
surface equation as $\beta=0$. In the case that $\beta\neq 0$, a
minimal surface with constant K\"ahler angle ( an infinitesimally
holomorphic immersion) satisfies the equation, especially two
kinds of important surfaces, i.e. holomorphic curves and special
Lagrangian surfaces satisfy the equation. However there are many
$\beta$-symplectic critical surfaces which are not minimal
surfaces.

Existence of holomorphic curves in a K\"ahler surface is a fundamental problem in differential geometry. It is known that a closed symplectic minimal surface in a compact K\"haler-Einstein surface with nonnegative scalar curvature is holomorphic (\cite{W1}). On the contrary, C. Arezzo (\cite{Arezzo}) constructed examples which shows that a strictly stable minimal surface in a K\"ahler-Einstein surface with negative scalar curvature may not be holomorphic.

Our goal in this paper is to start a program to deform $\beta$-symplectic critical surface from a minimal surface ($\beta=0$) to a holomorphic curve ($\beta=\infty$) by using continuity method. Actually, one can check that if we can deform the $\beta$-symplectic critical surface $\Sigma_{\beta}$ for $\beta\in [0,\infty)$ with uniformly bounded $L_{\beta}(\Sigma_{\beta})$ to a limit $\Sigma_{\infty}$ smoothly, then $\Sigma_{\infty}$ is a holomorphic curve in $M$, regardless of the sign of the scalar curvature of $M$. We do not need $M$ to be a K\"ahler-Einstein surface.

We first examine the properties of $\beta$-symplectic critical surfaces for each fixed $\beta$. We will see that $\beta$-symplectic critical surfaces share many properties with
minimal surfaces (c.f. \cite{CW}, \cite{HL}, \cite{MW},
\cite{W1},\cite{W2}). We derive an equation for the K\"ahler angle
of a $\beta$-symplectic critical surface in a K\"ahler-Einstein
surface.
\begin{theorem}\label{Lapan}
If $M$ is a K\"ahler-Einstein surface and $\Sigma$ is a
$\beta$-symplectic critical surface, then we have
\begin{eqnarray*}
\Delta\cos\alpha
&=&\frac{2\beta\sin^2\alpha}{\cos\alpha
(\cos^2\alpha+\beta\sin^2\alpha)}|\nabla\alpha|^2-2\cos\alpha|\nabla\alpha|^2-\frac{K}{4}\frac{\cos^3\alpha\sin^2\alpha}{\cos^2\alpha+\beta\sin^2\alpha},
\end{eqnarray*} where $K$ is the scalar curvature of $M$.
\end{theorem}
The theorem yields that a $\beta$-symplectic critical surface in a
K\"ahler-Einstein surface with nonnegative scalar curvature is
holomorphic.

\vspace{.1in}

Then we calculate the second variation formula of the functionals
$L_\beta$. As a corollary (c.f. \cite{HL}, \cite{MW}), we show that

\begin{theorem}
Let $M$ be a K\"ahler surface with positive scalar curvature $R$.
If $\Sigma$ is a stable $\beta$-symplectic critical surface in $M$
with $\beta\geq 0$ and $\chi(\nu)\geq g$, where $\chi(\nu)$ is the
Euler characteristic of the normal bundle $\nu$ of $\Sigma$ in $M$
and $g$ is the genus of $\Sigma$, then $\Sigma$ is a holomorphic
curve.
\end{theorem}

In order to proceed the continuity method, we define the set
\begin{equation*}
    S:=\{\beta\in[0,\infty)\mid \exists\ a\ strictly\ stable\ \beta-symplecitc\ critical\
    surface\}.
\end{equation*}
As a first attempt to the continuity method and as an application of the second variation formula, we prove that

\begin{theorem}
The set $S$ is open in $[0,\infty)$.
\end{theorem}

The existence part will depend on a theorem proved by B. White (\cite{White}) combining with Implicit Function Theorem. White's theorem tells us that the Jacobi operator is a Fredholm map with Fredholm index 0. The stability assumption implies that the Jacobi operator is injective, thus isomorphism. The stability of the solution follows from the continuity of the Jacobi operator with respect to $\beta$.
Of course, the closeness (or, the compactness) part is more dedicate. We will attach this problem in a subsequent paper.

As a verification of our idea, we study the rotationally symmetric graphic $\beta$-symplectic critical surface in ${\mathbb C}^2$ which is given by
\begin{equation*}
    F(r,\theta)=(r\cos\theta,r\sin\theta,f(r), g(r)).
\end{equation*}
We proved that such a surface is a $\beta$-symplectic critical surface if and only if
\begin{eqnarray*}\label{rot1}
 \left\{
  \begin{array}{cc}
   rf'(1+(f')^2+(g')^2)^{\frac{\beta-1}{2}} & = C_1, \\
     rg'(1+(f')^2+(g')^2)^{\frac{\beta-1}{2}} & = C_2.
  \end{array}
\right.
\end{eqnarray*}
By analyzing the ODE system carefully , we give the following asymptotic expansions:

\begin{theorem}
For any $\beta>0$ and $\varepsilon>0$, the equations
\begin{eqnarray*}
 \left\{
  \begin{array}{cc}
   rf'(1+(f')^2+(g')^2)^{\frac{\beta-1}{2}} & = 1, \\
     rg'(1+(f')^2+(g')^2)^{\frac{\beta-1}{2}} & = 1,\\
     f(\varepsilon) &=f_0,\\
     g(\varepsilon) &=g_0,
  \end{array}
\right.
\end{eqnarray*}
 have a unique $C^\infty$-solution on $[\varepsilon, +\infty)$. Moreover, $f'=g'$. As $r\to\infty$, we have the asymptotic expansion,
 $$f'=\frac{1}{r}-\frac{\beta-1}{r^3}+o(r^{-3}).$$
 As $r\to 0$, we have the asymptotic expansion,
 $$f'=2^{\frac{1-\beta}{2\beta}}r^{-\frac{1}{\beta}}-\frac{\beta-1}{\beta}2^{-\frac{3\beta+1}{2\beta}}r^{\frac{1}{\beta}}+o(r^{\frac{1}{\beta}}).$$
\end{theorem}

Furthermore, when $\beta\to0$, we show that $\Sigma_{\beta}$
converges to the catenoid locally, while when $\beta\to\infty$,
$\Sigma_{\beta}$ converges to a flat plane locally. This coincides
with our expectation.

The following sections are organized as follows: in Section 2, we
derive the Euler-Lagrange equation for the $L_{\beta}$ functional
and the elliptic equation satisfied by the K\"aher angle, we also
derive Webster's formula for $\beta$-symplectic critical surfaces;
in Section 3, we prove that the equation is an elliptic system
module tangential diffeomorphisms; in Section 4, we compute the
second variation for $L_{\beta}$ functional; in Section 5, we
prove the openness of the set of $\beta$ where there is a stable
$\beta$-symplectic critical surface; in the last section, we study
the rotationally symmetric $\beta$-symplectic critical surfaces in
${\mathbb C}^2$.

\subsection*{Acknowledgements} Partial of the work are carried out when the first and the third authors
are visiting the Abdus Salam International Centre for Theoretical Physics (ICTP). They would like to thank Professor Claudio Arezzo for invaluable discussions about this problem.
\vspace{.1in}

\section{The Euler-Lagrange Equation}

Assume that $\phi_t: \Sigma\to M$ is a one-parameter family of
immersions and $\frac{\partial\phi_t}{\partial
t}\left|_{t=0}\right.=\textbf{X}$, $\textbf{X}$ is the variational
vector field of $\Sigma$.  We denote by $\overline\nabla$ the
covariant derivative and by $K$ the Riemannian curvature tensor on
$M$. Furthermore, $\nabla, R$ denote the covariant derivative and
the Riemannian curvature tensor of the induced metric $g$ on the
surface $\Sigma$.
\begin{theorem}\label{th1}
Let $M$ be a K\"ahler surface. The first variational formula of
the functional $L_\beta$ is, for any smooth normal vector field
$\textbf{X}$ on $\Sigma$,
\begin{eqnarray}
\delta_X L_\beta &=& -(\beta+1)\int_{\Sigma}\frac{\textbf{X}\cdot
\textbf{H}}{\cos^\beta\alpha}
d\mu+\beta(\beta+1)\int_{\Sigma}\frac{\textbf{X}\cdot
(J(J\nabla\cos\alpha)^\top))^\bot}{\cos^{\beta+3}\alpha}d\mu,
\end{eqnarray}
where $\textbf{H}$ is the mean curvature vector of $\Sigma$ in
$M$, and $()^\top$ means tangential components of $()$, $()^\bot$
means the normal components of $()$. The Euler-Lagrange equation
of the functional $L_{\beta}$ for $\beta\neq -1$ is
\begin{equation}\label{betaequ}
\cos^3\alpha
\textbf{H}-\beta(J(J\nabla\cos\alpha)^\top)^\bot=0.
\end{equation}
\end{theorem}

\vspace{.1in}

\textbf{Proof:} Let $\{x_i\}$ be the local normal coordinates
around a fixed point $p$ on $\Sigma$. The induced metric on
$\phi_{t}(\Sigma)$ is
\begin{equation*}
    g_{ij}(t)=\langle\frac{\partial \phi_t}{\partial x_i},\frac{\partial \phi_t}{\partial x_j}\rangle.
\end{equation*}
For simplicity, we denote $\frac{\partial \phi_0}{\partial
x_i}$ by $e_i$ and $g_{ij}(t)$ by $g_{ij}$. A direct
calculation gives
\begin{equation}\label{e2.3}
    \frac{\partial}{\partial t}\mid _{t=0}g_{ij}=\langle \overline{\nabla}_{e_i}\textbf{X}, e_j\rangle
    +\langle e_i,\overline{\nabla}_{e_j}\textbf{X}\rangle.
\end{equation}
From the definition of K\"ahler angle we have
$$\cos\alpha_t=\frac{\omega(\partial\phi_t/\partial x^1,
\partial\phi_t/\partial x^2)}{\sqrt{\det(g_t)}},
$$ where $\det (g_t)$ is the determinant of the metric $(g_t)$.
Now the functional can be written as
$$L_\beta(\phi_t)=\int_{\Sigma}\frac{\det^{(\beta+1)/2}(g_t)}{\omega^\beta(\partial\phi_t/\partial x^1,
\partial\phi_t/\partial x^2)} dx^1\wedge dx^2.$$ Thus,
\begin{eqnarray}\label{e2.4}
 & &\frac{d}{dt}\left|_{t=0}\right.L_\beta(\phi_t)\nonumber\\
&=& (\beta+1)\int_{\Sigma}\frac{\langle
\overline{\nabla}_{e_i}\textbf{X},
e_i\rangle}{\cos^\beta\alpha}d\mu
-\beta\int_{\Sigma}\frac{\omega(\overline\nabla_{e_1}
\textbf{X}, e_2)+ \omega(e_1, \overline\nabla_{e_2}
\textbf{X})}{\cos^{\beta+1}\alpha}d\mu.
\end{eqnarray}
Now suppose $\textbf{X}$ is a normal vector field, then
(\ref{e2.4}) can be written as
\begin{eqnarray}\label{e2.5}
 & &\frac{d}{dt}\left|_{t=0}\right.L_\beta(\phi_t)\nonumber\\
 &=&\int_{\Sigma}\large
(\frac{\beta+1}{2}\frac{\partial_t g_{ij}|_{t=0}
g^{ij}}{\cos^\beta\alpha}-\beta\frac{\partial_t\omega(\partial\phi_t/\partial
x^1,
\partial\phi_t/\partial x^2)|_{t=0}}{\cos^{\beta+1}\alpha}\large )d\mu \nonumber\\
&=& -(\beta+1)\int_{\Sigma}\frac{\textbf{X}\cdot
\textbf{H}}{\cos^\beta\alpha}d\mu
-\beta\int_{\Sigma}\frac{\omega(\overline\nabla_{e_1}
\textbf{X}, e_2)+ \omega(e_1, \overline\nabla_{e_2}
\textbf{X})}{\cos^{\beta+1}\alpha}d\mu.
\end{eqnarray}
Since $\Sigma$ is closed, applying the Stokes formula, we obtain
\begin{eqnarray*}
& &-\beta\int_{\Sigma}\frac{\omega(\overline\nabla_{e_1}
\textbf{X}, e_2)+
\omega(e_1, \overline\nabla_{e_2} \textbf{X})}{\cos^{\beta+1}\alpha}d\mu\\
 &=&-\beta\int_{\Sigma}\frac{e_1(\omega(\textbf{X}, e_2))-\omega(\textbf{X},
\overline\nabla_{e_1}\overline\nabla_{e_2} F)}{\cos^{\beta+1}\alpha}d\mu\\
&&-\beta\int_{\Sigma}\frac{e_2(\omega(e_1,
\textbf{X}))-\omega(
\overline\nabla_{e_2}\overline\nabla_{e_1} F,
\textbf{X})}{\cos^{\beta+1}\alpha} d\mu\\
&=&-\beta\int_{\Sigma}\frac{\nabla_{e_1}(\omega(\textbf{X},
e_2))+\nabla_{e_2}(\omega(e_1, \textbf{X}))}{\cos^{\beta+1}\alpha}d\mu\\
&=&-\beta(\beta+1)\int_{\Sigma}\frac{\omega(\textbf{X},
e_2)\nabla_{e_1}\cos\alpha+\omega(e_1,
\textbf{X})\nabla_{e_2}\cos\alpha}{\cos^{\beta+2}\alpha} d\mu,
\end{eqnarray*} where we have used the fact that $\omega$ is
parallel. Since $\omega(\textbf{X}, e_2)=-\langle\textbf{X}, Je_2\rangle$,
$\omega(e_1, \textbf{X})=\langle\textbf{X}, Je_1\rangle$ and
\begin{eqnarray*}
(J\nabla\cos\alpha)^\top&=&(Je_1\nabla_{e_1}\cos\alpha+Je_2\nabla_{e_2}\cos\alpha)^\top\\
&=&\langle Je_1,e_2\rangle e_2\nabla_{e_1}\cos\alpha+\langle
Je_2,e_1\rangle e_1\nabla_{e_2}\cos\alpha\\
&=& (e_2\nabla_{e_1}\cos\alpha-
e_1\nabla_{e_2}\cos\alpha)\cos\alpha,
\end{eqnarray*} so
\begin{eqnarray*}
\omega(X, e_2)\nabla_{e_1}\cos\alpha+\omega(e_1,
X)\nabla_{e_2}\cos\alpha=-\frac{\textbf{X}\cdot
(J(J\nabla\cos\alpha)^\top))^\bot}{\cos\alpha}.
\end{eqnarray*} Therefore, we have
\begin{eqnarray*}
\frac{d}{dt}\left|_{t=0}\right.L_\beta(\phi_t)&=&-(\beta+1)\int_{\Sigma}\frac{\textbf{X}\cdot
\textbf{H}}{\cos^\beta\alpha}d\mu
\\&&+\beta(\beta+1)\int_\Sigma\frac{\textbf{X}\cdot
(J(J\nabla\cos\alpha)^\top))^\bot}{\cos^{\beta+3}\alpha} d\mu.
\end{eqnarray*}
\hfill Q. E. D.

\vspace{.1in}

We express $(J(J\nabla\cos\alpha)^\top))^\bot$ at a fixed point
$p$ in a local frame. If we assume the K\"ahler form is
self-dual, then $J$ has the form

\begin{eqnarray}\label{e15}
J=\left (\begin{array}{clcr} 0 &x &y &z \\
-x &0 &z &-y\\
-y &-z &0 &x\\
-z &y &-x &0 \end{array}\right),
\end{eqnarray} where $x^2+y^2+z^2=1.$
By the definition of the K\"ahler angle, we know that
$$x=\cos\alpha=\omega(e_1, e_2)=\langle Je_1, e_2\rangle.$$
Then,
\begin{eqnarray*}
(J(J\nabla\cos\alpha)^\top))^\bot &=&
(J(\cos\alpha\partial_1\cos\alpha
e_2-\cos\alpha\partial_2\cos\alpha e_1) )^\bot \\
&=&-\cos\alpha\sin\alpha\partial_1\alpha
(Je_2)^\bot+\cos\alpha\sin\alpha\partial_2\alpha (Je_1)^\bot\\
&=&-\cos\alpha\sin\alpha\partial_1\alpha(ze_3-ye_4)+\cos\alpha\sin\alpha\partial_2\alpha(ye_3+ze_4).
\end{eqnarray*}
Thus, by (\ref{betaequ}) we get that
\begin{eqnarray}\label{e2}
H^3&=&\beta\frac{\sin\alpha}{\cos^2\alpha}(y\partial_2\alpha-z\partial_1\alpha);\\
H^4&=&\beta\frac{\sin\alpha}{\cos^2\alpha}(y\partial_1\alpha+z\partial_2\alpha).
\end{eqnarray}

\begin{theorem}[(cf. \cite{HL})]
If $\Sigma$ is a closed symplectic surface which is smoothly
immersed in $M$ with the K\"ahler angle $\alpha$, then $\alpha$
satisfies the following equation ,
\begin{eqnarray}\label{cos2}
\Delta\cos\alpha &=& \cos\alpha(-|h^3_{1k}-h^4_{2k}|^2-
|h^4_{1k}+h^3_{2k}|^2)\nonumber
\\ &&+\sin\alpha(H^4_{,1}+H^3_{,2})-\frac{\sin^2\alpha}{\cos\alpha}
(K_{1212}+K_{1234}).
\end{eqnarray} where $K$ is the curvature operator of $M$ and
$H^\alpha_{,i}=\langle\overline\nabla_{e_i}^N {\textbf H}, e_\alpha\rangle$.
\end{theorem}

\begin{theorem}
Suppose that $M$ is K\"ahler surface and $\Sigma$ is a
$\beta$-symplectic critical surface in $M$ with K\"ahler angle
$\alpha$, then $\cos\alpha$ satisfies,
\begin{eqnarray}\label{e12}
\Delta\cos\alpha &=&\frac{2\beta\sin^2\alpha}{\cos\alpha
(\cos^2\alpha+\beta\sin^2\alpha)}|\nabla\alpha|^2-2\cos\alpha|\nabla\alpha|^2\nonumber\\
&&-\frac{\cos^2\alpha\sin^2\alpha}{\cos^2\alpha+\beta\sin^2\alpha}
Ric(Je_1, e_2).
\end{eqnarray}
\end{theorem}

{\it Proof.} We will compute pointwise. For a fixed point $p\in\Sigma$, we can choose the local frame such that at $p$, $y=\sin\alpha$ and $z=0$. For a $\beta$-symplectic critical surface $\Sigma$,
if we set $\textbf{V}=\nabla_{e_2}\alpha e_3+\nabla_{e_1}\alpha
e_4$, then we have
\begin{equation}\label{equation2}
    \textbf{H}=\beta\frac{\sin^2\alpha}{\cos^2\alpha}\textbf{V}.
\end{equation}
It is easy to check that (see (2.5) of \cite{HL})
\begin{equation}\label{e.alpha}
\nabla_{e_1}\alpha=-(h^4_{11}+h^3_{12}), \ \ \nabla_{e_2}\alpha=-(h^4_{12}+h^3_{22}).
\end{equation}
By direct computation, we have at $p$,
\begin{eqnarray*}
(h^3_{1k}-h^4_{2k})^2+(h^4_{1k}+h^3_{2k})^2 &=&
|\textbf{H}|^2+2|\textbf{V}|^2+2\textbf{H}\cdot \textbf{V} \\ &=&
(\beta^2\frac{\sin^4\alpha}{\cos^4\alpha}+2+2\beta\frac{\sin^2\alpha}{\cos^2\alpha})|\textbf{V}|^2
\\ &=& \frac{\beta^2\sin^4\alpha+2\cos^4\alpha+2\beta\sin^2\alpha\cos^2\alpha}{\cos^4\alpha} |\nabla\alpha|^2,
\end{eqnarray*}
and
\begin{eqnarray}\label{e5}
\sin\alpha(H^4_{, 1}+H^3_{,2})
&=&\sin\alpha(\langle\overline\nabla_{e_1}\textbf{H}, e_4\rangle+\langle\overline\nabla_{e_2}\textbf{H}, e_3\rangle)\nonumber\\
&=&\sin\alpha(\partial_1(H^4)+H^3\langle\overline\nabla_{e_1}e_3, e_4\rangle+\partial_2(H^3)+H^4\langle\overline\nabla_{e_2}e_4, e_3\rangle)\nonumber\\
&=&\beta\sin\alpha \{\partial_1[\frac{\sin\alpha}{\cos^2\alpha}(y\partial_1\alpha+z\partial_2\alpha)]+\partial_2[\frac{\sin\alpha}{\cos^2\alpha}(y\partial_2\alpha-z\partial_1\alpha)]\}\nonumber\\
&&+\sin\alpha(H^3\langle\overline\nabla_{e_1}e_3, e_4\rangle+H^4\langle\overline\nabla_{e_2}e_4, e_3\rangle)\nonumber\\
&=&\beta\sin\alpha[\partial_1(\frac{\sin\alpha}{\cos^2\alpha}\partial_1\alpha)y+\frac{\sin\alpha}{\cos^2\alpha}\partial_1\alpha\partial_1 y\nonumber\\
&&+\partial_1(\frac{\sin\alpha}{\cos^2\alpha}\partial_2\alpha)z+\frac{\sin\alpha}{\cos^2\alpha}\partial_2\alpha\partial_1 z\nonumber\\
&&+\partial_2(\frac{\sin\alpha}{\cos^2\alpha}\partial_2\alpha)y+\frac{\sin\alpha}{\cos^2\alpha}\partial_2\alpha\partial_2 y\nonumber\\
&&-\partial_2(\frac{\sin\alpha}{\cos^2\alpha}\partial_1\alpha)z-\frac{\sin\alpha}{\cos^2\alpha}\partial_1\alpha\partial_2 z]\nonumber\\
&&+\sin\alpha(H^3\langle\overline\nabla_{e_1}e_3, e_4\rangle+H^4\langle\overline\nabla_{e_2}e_4, e_3\rangle)\nonumber\\
&=&\beta\sin^2\alpha[\partial_1(\frac{\sin\alpha}{\cos^2\alpha}\partial_1\alpha)+\partial_2(\frac{\sin\alpha}{\cos^2\alpha}\partial_2\alpha)]\nonumber\\
&&+\beta\frac{\sin^2\alpha}{\cos^2\alpha}(\partial_1\alpha\partial_1y+\partial_2\alpha\partial_2y)
  +\beta\frac{\sin^2\alpha}{\cos^2\alpha}(\partial_2\alpha\partial_1 z-\partial_1\alpha\partial_2 z)\nonumber\\
&&+\sin\alpha(H^3\langle\overline\nabla_{e_1}e_3, e_4\rangle+H^4\langle\overline\nabla_{e_2}e_4, e_3\rangle)\nonumber\\
&=&\beta\frac{\sin^3\alpha}{\cos^2\alpha}\Delta\alpha+\beta\frac{\sin^2\alpha(1+\sin^2\alpha)}{\cos^3\alpha}|\nabla\alpha|^2\nonumber\\
&&+\beta\frac{\sin^2\alpha}{\cos^2\alpha}(\partial_1\alpha\partial_1y+\partial_2\alpha\partial_2y)
  +\beta\frac{\sin^2\alpha}{\cos^2\alpha}(\partial_2\alpha\partial_1 z-\partial_1\alpha\partial_2 z)\nonumber\\
&&+\beta\frac{\sin^3\alpha}{\cos^2\alpha}(\partial_2\alpha\langle\overline\nabla_{e_1}e_3, e_4\rangle+\partial_1\alpha\langle\overline\nabla_{e_2}e_4, e_3\rangle).
\end{eqnarray}

Now we begin to compute
$\partial_1\alpha\partial_1y+\partial_2\alpha\partial_2y$ and
$\partial_2\alpha\partial_1 z-\partial_1\alpha\partial_2 z$. Note
that $y=\langle Je_1, e_3\rangle$. Then
\begin{eqnarray*}
\partial_1\alpha\partial_1y+\partial_2\alpha\partial_2y
&=&\partial_1\alpha (\langle J\overline\nabla_{e_1}e_1,e_3\rangle+\langle Je_1,\overline\nabla_{e_1}e_3\rangle)\\
&&+\partial_2\alpha (\langle J\overline\nabla_{e_2}e_1,e_3\rangle+\langle Je_1,\overline\nabla_{e_2}e_3\rangle)\\
&=&\partial_1\alpha(-\cos\alpha h^4_{11}-\cos\alpha h^3_{12})
      +\partial_2\alpha (-\cos\alpha h^4_{12}-\cos\alpha h^3_{22})\\
&=&\cos\alpha(-h^4_{11}-h^3_{12})\partial_1\alpha+\cos\alpha (-h^4_{12}-h^3_{22})\partial_2\alpha \\
&=&\cos\alpha|\nabla\alpha|^2.
\end{eqnarray*}
Similarly, since $z=\langle Je_1, e_4\rangle$, we have
\begin{eqnarray*}
\partial_2\alpha\partial_1 z-\partial_1\alpha\partial_2 z
&=&\partial_2\alpha (\langle J\overline\nabla_{e_1}e_1,e_4\rangle+\langle Je_1,\overline\nabla_{e_1}e_4\rangle)\\
&&-\partial_1\alpha (\langle J\overline\nabla_{e_2}e_1,e_4\rangle+\langle Je_1,\overline\nabla_{e_2}e_4\rangle)\\
&=&\partial_2\alpha(-\cos\alpha h^4_{12}+\sin\alpha\langle\overline\nabla_{e_1}e_4, e_3\rangle+\cos\alpha h^3_{11})\\
&&-\partial_1\alpha(-\cos\alpha h^4_{22}+\sin\alpha\langle\overline\nabla_{e_2}e_4, e_3\rangle+\cos\alpha h^3_{12})\\
&=&\partial_2\alpha(H^3\cos\alpha+\cos\alpha\partial_2\alpha)+\partial_1\alpha(H^4\cos\alpha+\cos\alpha\partial_1\alpha)\\
&&+\sin\alpha(\partial_2\alpha\langle\overline\nabla_{e_1}e_4, e_3\rangle-\partial_1\alpha\langle\overline\nabla_{e_2}e_4, e_3\rangle)\\
&=&\frac{\cos^2\alpha+\beta\sin^2\alpha}{\cos\alpha}|\nabla\alpha|^2+\sin\alpha(\partial_2\alpha\langle\overline\nabla_{e_1}e_4,
e_3\rangle-\partial_1\alpha\langle\overline\nabla_{e_2}e_4,
e_3\rangle).
\end{eqnarray*}
Putting this equation into (\ref{e5}) we get that
\begin{equation}\label{e6}
   \sin\alpha(H^4_{, 1}+H^3_{,2})=-\beta\frac{\sin^2\alpha}{\cos^2\alpha}\Delta\cos\alpha
    +\beta\frac{\sin^2\alpha(2+\beta\sin^2\alpha)}{\cos^3\alpha}|\nabla\alpha|^2.
\end{equation}
Putting these two equations into (\ref{cos2}), we obtain that,
\begin{eqnarray*}
\Delta\cos\alpha &=&-
\frac{\beta^2\sin^4\alpha+2\cos^4\alpha+2\beta\sin^2\alpha\cos^2\alpha}{\cos^3\alpha}
|\nabla\alpha|^2-\beta\frac{\sin^2\alpha}{\cos^2\alpha}\Delta\cos\alpha\\
&&+\beta\frac{\sin^2\alpha(2+\beta\sin^2\alpha)}{\cos^3\alpha}|\nabla\alpha|^2
-\frac{\sin^2\alpha}{\cos\alpha}(K_{1212}+K_{1234})\\
&=&\frac{2\beta\sin^2\alpha-2\cos^2\alpha(\cos^2\alpha+\beta\sin^2\alpha)}{\cos^3\alpha}
|\nabla\alpha|^2\\&&-\beta\frac{\sin^2\alpha}{\cos^2\alpha}\Delta\cos\alpha-\frac{\sin^2\alpha}{\cos\alpha}(K_{1212}+K_{1234}).
\end{eqnarray*} Therefore,

\begin{eqnarray*}
\Delta\cos\alpha &=&\frac{2\beta\sin^2\alpha}{\cos\alpha
(\cos^2\alpha+\beta\sin^2\alpha)}|\nabla\alpha|^2-2\cos\alpha|\nabla\alpha|^2\\
&&-\frac{\cos\alpha\sin^2\alpha}{\cos^2\alpha+\beta\sin^2\alpha}
(K_{1212}+K_{1234})\\
&=&\frac{2\beta\sin^2\alpha}{\cos\alpha
(\cos^2\alpha+\beta\sin^2\alpha)}|\nabla\alpha|^2-2\cos\alpha|\nabla\alpha|^2\\
&&-\frac{\cos^2\alpha\sin^2\alpha}{\cos^2\alpha+\beta\sin^2\alpha}
Ric(Je_1, e_2).
\end{eqnarray*}
The last equality used the fact that $K_{1212}+K_{1234}=\cos\alpha Ric(Je_1,e_2)$ (Lemma 3.2 in \cite{HL}).
\hfill Q. E. D.

\begin{corollary}\label{cor.angle}
Assume $M$ is K\"ahler-Einstein surface with scalar curvature $K$,
then $\cos\alpha$ satisfies,
\begin{equation*}\label{e14}
\Delta\cos\alpha =\frac{2\beta\sin^2\alpha}{\cos\alpha
(\cos^2\alpha+\beta\sin^2\alpha)}|\nabla\alpha|^2-2\cos\alpha|\nabla\alpha|^2
-\frac{K}{4}\frac{\cos^3\alpha\sin^2\alpha}{\cos^2\alpha+\beta\sin^2\alpha}.
\end{equation*}
\end{corollary}

\begin{corollary}
Any $\beta$-symplectic critical surface in a K\"ahler-Einstein
surface with nonnegative scalar curvature is a holomorphic curve for
$\beta\geq0$.
\end{corollary}

\vspace{.1in}

Similar to minimal surfaces (\cite{W1}, \cite{W2}) and symplectic
critical surfaces (\cite{HL}), a nonholomorphic $\beta$-symplectic
critical surface in a K\"ahler surface has at most finite complex
points, and the following formula can be derived.

\begin{theorem}\label{toplocal}
Suppose that $\Sigma$ is a non holomorphic $\beta$-symplectic
critical surface in a K\"ahler surface $M$. Then
$$
\chi (\Sigma)+\chi
(\nu)=-P,
$$
and
$$
c_1(M)([\Sigma])=-P,
$$
where $\chi(\Sigma)$ is the Euler characteristic of $\Sigma$,
$\chi(\nu)$ is the Euler characteristic of the normal bundle of
$\Sigma$ in $M$, $c_1(M)$ is the first Chern class of $M$,
$[\Sigma]\in H_2(M,{\bf Z})$ is the homology class of $\Sigma$ in
$M$, and $P$ is the number of complex tangent points.
\end{theorem}

\textbf{Proof :} Similar to the analysis in \cite{HL}, we know that the complex points are isolated. Set $$
g(\alpha)=\ln(\sin^2\alpha).$$ Then using the equation
(\ref{e12}), we obtain
\begin{eqnarray}\label{e13}
\Delta g(\alpha)
&=&-2|\nabla\alpha|^2-2\frac{\cos\alpha}{\sin^2\alpha}\Delta\cos\alpha
-4\frac{\cos^2\alpha}{\sin^2\alpha}|\nabla\alpha|^2 \nonumber\\
&=&-2|\nabla\alpha|^2-\frac{4\beta}{\cos^2\alpha+\beta\sin^2\alpha}|\nabla\alpha|^2\nonumber\\
&&+\frac{2\cos^2\alpha}{\cos^2\alpha+\beta\sin^2\alpha}(K_{1212}+K_{1234}).
\end{eqnarray} This equation is valid away from the complex
tangent points of $M$. By the Gauss equation and Ricci equation,
we have,
\begin{eqnarray*}
R_{1212}&=&K_{1212}+h^\alpha_{11}h^\alpha_{22}-(h^\alpha_{12})^2\\
R_{1234}&=&K_{1234}+h^3_{1k}h^4_{2k}-h^4_{1k}h^3_{2k},
\end{eqnarray*} where $R_{1212}$ is the curvature of $T\Sigma$ and $R_{1234}$ is the
curvature of the normal bundle $\nu$. Adding these two equations together, we
get that,
\begin{eqnarray*}
K_{1212}+K_{1234} &=& R_{1212}+R_{1234}-\frac{1}{2}|\textbf{H}|^2 \\
&&+\frac{1}{2}((h^3_{1k}-h^4_{2k})^2+(h^4_{1k}+h^3_{2k})^2)\\ &=&
R_{1212}+R_{1234}+|\textbf{V}|^2+\textbf{H}\cdot \textbf{V} \\ &=&
R_{1212}+R_{1234}+\frac{\cos^2\alpha+\beta\sin^2\alpha}{\cos^2\alpha}|\nabla
\alpha|^2.
\end{eqnarray*} Thus,
\begin{eqnarray*}
R_{1212}+R_{1234}&=&\frac{\cos^2\alpha+\beta\sin^2\alpha}{2\cos^2\alpha}\Delta
g(\alpha)+\frac{2\beta}{\cos^2\alpha}|\nabla\alpha|^2.
\end{eqnarray*}
Integrating the above equality over $\Sigma$, arguing as that in \cite{W1}, we
can obtain,
\begin{equation}\label{A}
    2\pi(\chi(T\Sigma)+\chi(\nu))=-2\pi P,
\end{equation}
where $\chi(T\Sigma)$ is the Euler characteristic of $\Sigma$ and
$\chi(\nu)$ is the Euler characteristic of the normal bundle
of $\Sigma$ in $M$, $P$ is the sum of the orders of complex
tangent points.

By (\ref{e13}) we also get that,
\begin{eqnarray*}
\Delta g(\alpha)&=& -2|\nabla\alpha|^2-\frac{4\beta}{\cos^2\alpha+\beta\sin^2\alpha}|\nabla\alpha|^2\nonumber\\
&&+\frac{2\cos^3\alpha}{\cos^2\alpha+\beta\sin^2\alpha}Ric(Je_1,
e_2).
\end{eqnarray*}
Note that $Ric(Je_1, e_2) d\mu_\Sigma$ is the pulled back Ricci 2-form of $M$ by the
immersion $F$ to $\Sigma$, i.e,
$$ F^\ast (Ric^M)=Ric(Je_1, e_2) d\mu_\Sigma.$$ Thus,
$$F^\ast (Ric^M)=(\frac{\cos^2\alpha+\beta\sin^2\alpha}{2\cos^3\alpha}\Delta g(\alpha)
+\frac{2\beta+\cos^2\alpha+\beta\sin^2\alpha}{\cos^3\alpha}|\nabla\alpha|^2)d\mu_\Sigma.
$$ Integrating it over $\Sigma$, we obtain that,
\begin{equation}\label{B}
    2\pi F^\ast c_1(M)[\Sigma]=-2\pi P.
\end{equation}
\hfill Q.E.D.

\begin{remark}
In \cite{HL} there was an error in the computations which was corrected above.
\end{remark}

\begin{remark}
The formulae (\ref{A}) and (\ref{B}) are usually called Webter's formula, and proved by Webster in \cite{Web} and \cite{Web2}.
\end{remark}

\vspace{.1in}

\section{Principal Symbol of the Equation (\ref{betaequ})}

In this section, we examine the principle symbol of the equation (\ref{betaequ}).
First note that, the equation (\ref{betaequ}) can be rewritten as
\begin{equation}\label{equation}
    \cos^2\alpha\textbf{H}=\beta(J(\nabla_{e_1}\cos\alpha e_2-\nabla_{e_2}\cos\alpha e_1))^{\perp}.
\end{equation}
For simplicity, we suppose that $\Sigma$ is a surface in
${\mathbb C}^{2}$. The general case is similar. In local
coordinate, we can express the surface as
\begin{eqnarray*}
  F: \Sigma & \longrightarrow & {\mathbb C}^{2}={\mathbb R}^{4} \\
   (x_1,x_2)&\longmapsto &   F(x_1,x_2)=(F^1(x_1,x_2),\cdots,F^{4}(x_1,x_2)).
\end{eqnarray*}
We will use the following conventions:
\begin{equation*}
    1\leq i,j,\cdots \leq 2, \ \ 3\leq \alpha,\beta,\cdots  \leq 4, \ \ 1\leq A,B,\cdots \leq 4.
\end{equation*}
The tangent space of $\Sigma$ at a fixed point $x\in\Sigma$ is
spanned by $\{e_1,e_2\}$ given by
\begin{equation}\label{eA.5}
    e_1=\frac{\partial F}{\partial x_1}=\frac{\partial F^A}{\partial x_1}E_A, \ \ \
     e_2=\frac{\partial F}{\partial x_2}=\frac{\partial F^A}{\partial x_2}E_A,
\end{equation}
where $\{E_1,\cdots,E_{4}\}$ is the standard orthonormal basis of
${\mathbb R}^{4}$. Therefore, the induced metric on $\Sigma$ is
given by
\begin{equation}\label{eA.6}
    g_{ij}=\langle e_i,e_j\rangle=\frac{\partial F^A}{\partial x_i}\frac{\partial F^A}{\partial x_j}.
\end{equation}
We can take the coordinates so that at the fixed point
$x\in\Sigma$, we have $g_{ij}(x)=\delta_{ij}$. We will also take
the standard complex structure $J$ on ${\mathbb C}^{2}$ given by
\begin{equation}\label{eA.7}
J=\left(
  \begin{array}{cccc}
    0 &  -1  &  0 & 0\\
    1  & 0 &  0 &  0 \\
     0 &  0  & 0 & -1 \\
    0  &   0 & 1 & 0 \\
  \end{array}
\right).
\end{equation}
Then we have
\begin{equation}\label{eA.8}
    \begin{cases}
        JE_{2k-1}=E_{2k},    \\
        JE_{2k}=-E_{2k-1}.
    \end{cases}
\end{equation}
Furthermore, we choose any orthonormal basis
$\{e_{\alpha}\}_{\alpha=3}^{4}$ of the normal space. Denote
\begin{equation}\label{eA.9}
    P=\cos^2\alpha\textbf{H}-\beta(J(\nabla_{e_1}\cos\alpha e_2-\nabla_{e_2}\cos\alpha e_1))^{\perp}.
\end{equation}
We will compute the principal symbol of $P$.

\vspace{.1in}

\noindent First we consider the principal part of $\textbf{H}$.
Note that by (\ref{eA.6}), we can easily see that the Christoffel
symbol of the induced metric is:
\begin{equation*}
    \Gamma_{ij}^{k}=\frac{1}{2}g^{kl}\left\{\frac{\partial g_{il}}{\partial x_j}+\frac{\partial g_{jl}}{\partial x_i}
           -\frac{\partial g_{ij}}{\partial x_l}\right\}
           =g^{kl}\frac{\partial^2 F^{B}}{\partial x_i\partial x_j}\frac{\partial F^{B}}{\partial x_l}.
\end{equation*}
Therefore, we have
\begin{eqnarray}\label{eA.10}
\textbf{H}=\Delta_{\Sigma}F
   &=&g^{ij}\left(\frac{\partial^2F}{\partial x_i\partial x_j}-\Gamma^{k}_{ij}\frac{\partial F}{\partial x_k}\right)\nonumber\\
   &=& g^{ij}\left(\frac{\partial^2F}{\partial x_i\partial x_j}
           -g^{kl}\frac{\partial^2 F^{B}}{\partial x_i\partial x_j}\frac{\partial F^{B}}{\partial x_l}\frac{\partial F}{\partial x_k}\right)\nonumber\\
   &=& g^{ij}\left(\frac{\partial^2F^A}{\partial x_i\partial x_j}
           -g^{kl}\frac{\partial F^A}{\partial x_k}\frac{\partial F^{B}}{\partial x_l}\frac{\partial^2 F^{B}}{\partial x_i\partial x_j}\right)E_A.
\end{eqnarray}
The linearization of the operator at $F$ in the direction $G$ is:
\begin{equation}\label{eA.11}
    D(\textbf{H})(F)G=g^{ij}\left(\frac{\partial^2 G^{A}}{\partial x_i\partial x_j}-
          g^{kl}\frac{\partial F^{A}}{\partial x_k}\frac{\partial F^{B}}{\partial x_l}\frac{\partial^2 G^{B}}{\partial x_i\partial x_j} \right)E_A
             +first\ order\  terms.
\end{equation}

\vspace{.1in}

\noindent Next, we will consider the second part of $P$. By
definition,
\begin{equation}\label{eA.12}
    \cos\alpha=\frac{\omega(e_1,e_2)}{\sqrt{det (g_{ij})}}=\frac{\langle Je_1,e_2\rangle}{\sqrt{det (g_{ij})}}.
\end{equation}
By (\ref{eA.5}) and (\ref{eA.7}), we have:
\begin{equation*}\label{eA.13}
    \langle Je_1,e_2\rangle=\sum_{k=1}^{2}\left(\frac{\partial F^{2k-1}}{\partial x_1}\frac{\partial F^{2k}}{\partial x_2}
       -\frac{\partial F^{2k}}{\partial x_1}\frac{\partial F^{2k-1}}{\partial x_2}\right).
\end{equation*}
Therefore, we have
\begin{equation*}
    \cos\alpha=\frac{\sum_{k=1}^{2}\left(\frac{\partial F^{2k-1}}{\partial x_1}\frac{\partial F^{2k}}{\partial x_2}
       -\frac{\partial F^{2k}}{\partial x_1}\frac{\partial F^{2k-1}}{\partial x_2}\right)}{\sqrt{det (g_{ij})}},
\end{equation*}
\begin{eqnarray*}\label{eA.14}
   & & \frac{\partial \cos\alpha}{\partial x_1} = \nonumber\\
   &=& \frac{1}{\sqrt{det (g_{ij})}}\left\{\sum_{k=1}^{2}\left(\frac{\partial^2 F^{2k-1}}{\partial x_1^2}\frac{\partial F^{2k}}{\partial x_2}
       +\frac{\partial F^{2k-1}}{\partial x_1}\frac{\partial^2 F^{2k}}{\partial x_1\partial x_2}
       -\frac{\partial^2 F^{2k}}{\partial x_1^2}\frac{\partial F^{2k-1}}{\partial x_2}
       -\frac{\partial F^{2k}}{\partial x_1}\frac{\partial^2 F^{2k-1}}{\partial x_1\partial x_2}\right)\right.\nonumber\\
   & &\ \ \ \ \ \ \ \ \ \ \ \ \ \  \left. -\frac{1}{2}\langle Je_1,e_2\rangle g^{ij}
         \left(\frac{\partial^2 F^A}{\partial x_1\partial x_i}\frac{\partial F^A}{\partial x_j}
            +\frac{\partial F^A}{\partial x_i}\frac{\partial^2 F^A}{\partial x_1\partial x_j}\right)\right\},\nonumber\\
\end{eqnarray*}
and
\begin{eqnarray*}\label{eA.15}
   & & \frac{\partial \cos\alpha}{\partial x_2} =\nonumber\\
   &=& \frac{1}{\sqrt{det (g_{ij})}}\left\{\sum_{k=1}^{2}\left(\frac{\partial^2 F^{2k-1}}{\partial x_1\partial x_2}\frac{\partial F^{2k}}{\partial x_2}
       +\frac{\partial F^{2k-1}}{\partial x_1}\frac{\partial^2 F^{2k}}{\partial x_2^2}
       -\frac{\partial^2 F^{2k}}{\partial x_1\partial x_2}\frac{\partial F^{2k-1}}{\partial x_2}
       -\frac{\partial F^{2k}}{\partial x_1}\frac{\partial^2 F^{2k-1}}{\partial x_2^2}\right)\right.\nonumber\\
   & &\ \ \ \ \ \ \ \ \ \ \ \ \ \  \left. -\frac{1}{2}\langle Je_1,e_2\rangle g^{ij}
         \left(\frac{\partial^2 F^A}{\partial x_2\partial x_i}\frac{\partial F^A}{\partial x_j}
            +\frac{\partial F^A}{\partial x_i}\frac{\partial^2 F^A}{\partial x_2\partial x_j}\right)\right\}.\nonumber\\
\end{eqnarray*}
By our choice of the frame, at the fixed point $x$, we have
\begin{eqnarray}\label{eA.16}
   & & (J(\nabla_{e_1}\cos\alpha e_2-\nabla_{e_2}\cos\alpha e_1))^{\perp}\nonumber\\
   &=& \frac{\partial \cos\alpha}{\partial x_1}(Je_2)^{\perp}-\frac{\partial \cos\alpha}{\partial x_2}(Je_1)^{\perp}\nonumber\\
   &=& \frac{\partial \cos\alpha}{\partial x_1}\langle Je_2,e_{\alpha}\rangle e_{\alpha}
           -\frac{\partial \cos\alpha}{\partial x_2}\langle Je_1,e_{\alpha}\rangle e_{\alpha}.
\end{eqnarray}

\vspace{.1in}

\noindent Notice that $\cos\alpha$, $g_{ij}$, $e_i$ and
$e_{\alpha}$ only involve first order derivatives of the immersion
$F$. Therefore, by (\ref{eA.11}) and (\ref{eA.16}), we know that
the linearization of the operator $P$ at $F$ in the direction $G$
(computed at the point $x$) is:
\begin{eqnarray}\label{eA.17}
   & &D(\textbf{P})(F)G \nonumber\\
   &=& \cos^2\alpha g^{ij}\left(\frac{\partial^2 G^{A}}{\partial x_i\partial x_j}-
          g^{kl}\frac{\partial F^{A}}{\partial x_k}\frac{\partial F^{B}}{\partial x_l}\frac{\partial^2 G^{B}}{\partial x_i\partial x_j} \right)E_A\nonumber\\
   & & -\beta\left\{\sum_{k=1}^{n}\left(\frac{\partial^2 G^{2k-1}}{\partial x_1^2}\frac{\partial F^{2k}}{\partial x_2}
       +\frac{\partial F^{2k-1}}{\partial x_1}\frac{\partial^2 G^{2k}}{\partial x_1\partial x_2}
       -\frac{\partial^2 G^{2k}}{\partial x_1^2}\frac{\partial F^{2k-1}}{\partial x_2}
       -\frac{\partial F^{2k}}{\partial x_1}\frac{\partial^2 G^{2k-1}}{\partial x_1\partial x_2}\right)\right.\nonumber\\
   & &\ \ \ \ \ \ \ \ \ \ \ \ \ \  \left. -\frac{1}{2}\langle Je_1,e_2\rangle g^{ij}
         \left(\frac{\partial^2 G^A}{\partial x_1\partial x_i}\frac{\partial F^A}{\partial x_j}
            +\frac{\partial F^A}{\partial x_i}\frac{\partial^2 G^A}{\partial x_1\partial x_j}\right)\right\}\langle Je_2,e_{\alpha}\rangle e_{\alpha}\nonumber\\
   & & +\beta\left\{\sum_{k=1}^{n}\left(\frac{\partial^2 G^{2k-1}}{\partial x_1\partial x_2}\frac{\partial F^{2k}}{\partial x_2}
       +\frac{\partial F^{2k-1}}{\partial x_1}\frac{\partial^2 G^{2k}}{\partial x_2^2}
       -\frac{\partial^2 G^{2k}}{\partial x_1\partial x_2}\frac{\partial F^{2k-1}}{\partial x_2}
       -\frac{\partial F^{2k}}{\partial x_1}\frac{\partial^2 G^{2k-1}}{\partial x_2^2}\right)\right.\nonumber\\
   & &\ \ \ \ \ \ \ \ \ \ \ \ \ \  \left. -\frac{1}{2}\langle Je_1,e_2\rangle g^{ij}
         \left(\frac{\partial^2 G^A}{\partial x_2\partial x_i}\frac{\partial F^A}{\partial x_j}
            +\frac{\partial F^A}{\partial x_i}\frac{\partial^2 G^A}{\partial x_2\partial x_j}\right)\right\}\langle Je_1,e_{\alpha}\rangle e_{\alpha}
           \nonumber\\
   & & + \ first\ order\  terms.
\end{eqnarray}
\noindent We will denote $G^T$ and $G^{\perp}$ the projection of
$G\in {\mathbb R}^{4}$ on the tangent bundle and normal bundle of
$\Sigma$ respectively. It is easy to see that
\begin{equation*}
    |G^{T}|^2=g^{kl}\langle G,\frac{\partial F}{\partial x_k}\rangle\langle G,\frac{\partial F}{\partial x_l}\rangle
         =g^{kl}G^{A}G^{B}\frac{\partial F^{A}}{\partial x_k}\frac{\partial F^{B}}{\partial x_l}.
\end{equation*}
Then we see that the principal symbol of $P$ is given by:
\begin{eqnarray}\label{eA.18}
   & & \langle\sigma(D(P))(x,\xi)G,G \rangle \nonumber\\
   &=& \cos^2\alpha g^{ij}\left(\xi_i\xi_j|G|^2-
          g^{kl}\frac{\partial F^{A}}{\partial x_k}\frac{\partial F^{B}}{\partial x_l}\xi_i\xi_jG^AG^B \right)\nonumber\\
   & & -\beta\left\{\sum_{k=1}^{n}\left(G^{2k-1}\frac{\partial F^{2k}}{\partial x_2}-G^{2k}\frac{\partial F^{2k-1}}{\partial x_2}\right)
            \langle Je_2,G^{\perp}\rangle \xi_1^2\right.\nonumber\\
   & & \left.+\sum_{k=1}^{n}\left(G^{2k-1}\frac{\partial F^{2k}}{\partial x_1}-G^{2k}\frac{\partial F^{2k-1}}{\partial x_1}\right)
            \langle Je_1,G^{\perp}\rangle\xi_2^2\right.\nonumber\\
   & & \left.-\sum_{k=1}^{n}\left[\left(G^{2k-1}\frac{\partial F^{2k}}{\partial x_1}-G^{2k}\frac{\partial F^{2k-1}}{\partial x_1}\right)
            \langle Je_2,G^{\perp}\rangle \right.\right.\nonumber\\
   & & \left.\left.  +\left(G^{2k-1}\frac{\partial F^{2k}}{\partial x_2}-G^{2k}\frac{\partial F^{2k-1}}{\partial x_2}\right)
            \langle Je_1,G^{\perp}\rangle\right]\xi_1\xi_2 \right.\nonumber\\
   & &  \left.+\cos\alpha g^{ij}G^A\left(\frac{\partial F^A}{\partial x_j}\langle Je_1,G^{\perp}\rangle \xi_2\xi_i
               -\frac{\partial F^A}{\partial x_j}\langle Je_2,G^{\perp}\rangle \xi_1\xi_i\right) \right\}.
\end{eqnarray}
By (\ref{eA.8}), we have
\begin{eqnarray}\label{eA.19}
  & & \langle\sigma(D(P))(x,\xi)G,G \rangle \nonumber\\
  &=& \cos^2\alpha|\xi|^2|G^{\perp}|^2\nonumber\\
     & & -\beta\left\{\left(-\langle Je_2,G^{\perp}\rangle\langle Je_2,G\rangle
         -\langle Je_1,e_2\rangle\langle G,e_1\rangle\langle Je_2,G^{\perp}\rangle\right) \xi_1^2\right.\nonumber\\
   & &  \left.+\left(-\langle Je_1,G^{\perp}\rangle\langle Je_1,G\rangle
         +\langle Je_1,e_2\rangle\langle G,e_2\rangle\langle Je_1,G^{\perp}\rangle\right) \xi_2^2\right.\nonumber\\
   & &  \left.+\left(\langle Je_2,G^{\perp}\rangle\langle Je_1,G\rangle+\langle Je_1,G^{\perp}\rangle\langle Je_2,G\rangle \right.\right.\nonumber\\
   & &   \left.\left.  +\langle Je_1,e_2\rangle\langle G,e_1\rangle\langle Je_1,G^{\perp}\rangle
            -\langle Je_1,e_2\rangle\langle G,e_2\rangle\langle Je_2,G^{\perp}\rangle\right)\xi_1\xi_2\right\}.
\end{eqnarray}
Note that $(Je_1)^T=\langle Je_1,e_2\rangle e_2$ and
$(Je_2)^T=-\langle Je_1,e_2\rangle e_1$. Thus we have
\begin{eqnarray}\label{eA.20}
   \langle\sigma(D(P))(x,\xi)G,G \rangle &=& \cos^2\alpha|\xi|^2|G^{\perp}|^2\nonumber\\
   & & +\beta\left(\langle G^{\perp},Je_2\rangle^2\xi_1^2-2\langle G^{\perp},Je_1\rangle\langle G^{\perp},Je_2\rangle\xi_1\xi_2
           +\langle G^{\perp},Je_1\rangle^2\xi_2^2\right)\nonumber\\
   &=& \left(\cos^2\alpha|G^{\perp}|^2+\beta\langle G^{\perp},Je_2\rangle^2\right)\xi_1^2
       -2\beta\langle G^{\perp},Je_1\rangle\langle G^{\perp},Je_2\rangle\xi_1\xi_2\nonumber\\
   & & +\left(\cos^2\alpha|G^{\perp}|^2+\beta\langle G^{\perp},Je_1\rangle^2\right)\xi_2^2.
\end{eqnarray}
The coefficient matrix is given by
\begin{equation}\label{eA.21}
  O=  \left(
  \begin{array}{cc}
    \cos^2\alpha|G^{\perp}|^2+\beta\langle G^{\perp},Je_2\rangle^2 & -\beta\langle G^{\perp},Je_1\rangle\langle G^{\perp},Je_2\rangle \\
    -\beta\langle G^{\perp},Je_1\rangle\langle G^{\perp},Je_2\rangle  & \cos^2\alpha|G^{\perp}|^2+\beta\langle G^{\perp},Je_1\rangle^2 \\
  \end{array}
\right).
\end{equation}
We have
\begin{eqnarray}\label{eA.22}
 det O
   &=& \cos^4\alpha|G^{\perp}|^4+\beta\cos^2\alpha|G^{\perp}|^2(\langle G^{\perp},Je_1\rangle^2+\langle G^{\perp},Je_2\rangle^2)\nonumber\\
   &=& \cos^2\alpha|G^{\perp}|^2\left[\cos^2\alpha|G^{\perp}|^2+\beta(\langle G^{\perp},Je_1\rangle^2+\langle G^{\perp},Je_2\rangle^2)\right].
\end{eqnarray}

\vspace{.1in}

From (\ref{eA.22}), we see that \emph{if $\beta\geq0$, then $det O\geq0$, with equality holds if
and only if $\cos^2\alpha|G^{\perp}|^2=0$. } In particular,

\begin{proposition}\label{prop3.1}
For a $\beta$-symplectic critical surface with $\beta\geq0$, the
Euler-Lagrange equation (\ref{betaequ}) is an elliptic system modulo
tangential diffeomorphisms of $\Sigma$.
\end{proposition}

\vspace{.1in}

\section{Second Variation Formula}

In this section, we will compute the second variation of the
functional $L_{\beta}$ for a $\beta$-symplectic critical surface.
To this end, we use the notations as in Section 2. Recall that
\begin{equation*}
    L_\beta(\phi)=\int_\Sigma\frac{1}{\cos^\beta\alpha}d\mu,
\end{equation*}
where $\phi:\Sigma\to M$ is a smooth symplectic immersion. Let
$\phi(\Sigma)$ be a $\beta$-symplectic critical surface. We
consider a smooth family of symplectic immersions
\begin{equation*}
    \phi_{t,\epsilon}:\Sigma\times(-\delta,\delta)\times (-a,a)\to M
\end{equation*}
with $\phi_{0,0}=\phi$. Since symplectic is an open condition, we
can assume that $\phi_{t,\epsilon}(\Sigma)$ is symplectic property for
every $t\in(-\delta,\delta)$ and every $\epsilon\in(-a,a)$, that
means that $\phi_{t,\epsilon}$ is a symplectic variation of
$\phi$. For symplicity, we write
$\phi_{0,0}(\Sigma)=\phi(\Sigma)=\Sigma$. Let
\begin{equation*}
    \frac{\partial \phi_{t,0}}{\partial t}\mid_{t=0}=\textbf{X}, \ \ \frac{\partial \phi_{0,\epsilon}}{\partial \epsilon}\mid_{\epsilon=0}=\textbf{Y},
    \ \ \mbox{and} \ \ \frac{\partial^2 \phi_{t,\epsilon}}{\partial t\partial \epsilon}\mid_{t=0,\epsilon=0}=\textbf{Z}.
\end{equation*}

\begin{lemma}\label{lem4.1}
We have the general second variation formula:
\begin{eqnarray}\label{2nd1}
   & &\frac{\partial^2}{\partial t\partial\epsilon}\mid _{t=0,\epsilon=0} L_{\beta}(\phi_{t,\epsilon})\nonumber\\
   &=& (\beta+1)\int_{\Sigma}\frac{\langle\overline{\nabla}_{e_i}^{\perp}\textbf{X},\overline{\nabla}_{e_i}^{\perp}\textbf{Y}\rangle}
               {\cos^{\beta}\alpha} d\mu
         -(\beta+1)\int_{\Sigma}\frac{K(\textbf{X},e_i,\textbf{Y},e_i)}{\cos^{\beta}\alpha}d\mu\nonumber\\
   & &   +(\beta+1)^2\int_{\Sigma}\frac{div_{\Sigma}\textbf{X}div_{\Sigma}\textbf{Y}}{\cos^{\beta}\alpha} d\mu\nonumber\\
   & & -\frac{\beta+1}{2}\int_{\Sigma}\frac{\langle e_i,\overline\nabla_{e_j}\textbf{X}\rangle\langle e_j,\overline\nabla_{e_i}\textbf{Y}\rangle+\langle e_j,\overline\nabla_{e_i}\textbf{X}\rangle\langle e_i,\overline\nabla_{e_j}\textbf{Y}\rangle}{\cos^{\beta}\alpha} d\mu\nonumber\\
   & & -\beta(\beta+1)\int_{\Sigma}\frac{div_{\Sigma}\textbf{X}(\omega(\overline{\nabla}_{e_1}\textbf{Y}, e_2)+\omega(e_1,\overline{\nabla}_{e_2}\textbf{Y}))}{\cos^{\beta+1}\alpha}d\mu\nonumber\\
   & & -\beta(\beta+1)\int_{\Sigma}\frac{div_{\Sigma}\textbf{Y}(\omega(\overline{\nabla}_{e_1}\textbf{X}, e_2)+\omega(e_1,\overline{\nabla}_{e_2}\textbf{X}))}{\cos^{\beta+1}\alpha}d\mu\nonumber\\
   & &  -\beta\int_{\Sigma}\frac{\omega(K(\textbf{Y},e_1)\textbf{X},e_2)
   \omega(e_1,K(\textbf{Y},e_2)\textbf{X})}
                {\cos^{\beta+1}\alpha} d\mu\nonumber\\
   & &  -\beta\int_{\Sigma}\frac{\omega(\overline{\nabla}_{e_1}\textbf{X},\overline{\nabla}_{e_2}\textbf{Y})
             +\omega(\overline{\nabla}_{e_1}\textbf{Y},\overline{\nabla}_{e_2}\textbf{X})}
                {\cos^{\beta+1}\alpha} d\mu\\
   & & +\beta(\beta+1)\int_{\Sigma}\frac{(\omega(\overline{\nabla}_{e_1}\textbf{X}, e_2)
          +\omega(e_1,\overline{\nabla}_{e_2}\textbf{X}))(\omega(\overline{\nabla}_{e_1}\textbf{Y}, e_2)
          +\omega(e_1,\overline{\nabla}_{e_2}\textbf{Y}))}{\cos^{\beta+2}\alpha}d\mu\nonumber,
\end{eqnarray}
where $\{e_1,e_2\}$ is a local orthonormal basis of $\Sigma$,
$div_{\Sigma}\textbf{X}=\langle \overline{\nabla}_{e_i}\textbf{X},
e_i\rangle$ and $K$ is the curvature tensor on $M$.
\end{lemma}

\vspace{.1in}

\textbf{Proof:} As in Section 2, let $\{x_i\}$ be the local normal
coordinates around a fixed point $p$ on $\Sigma$. The induced
metric on $\phi_{t,\epsilon}(\Sigma)$ is
\begin{equation*}
    g_{ij}(t,\epsilon)=\langle\frac{\partial \phi_{t,\epsilon}}{\partial x_i},\frac{\partial \phi_{t,\epsilon}}{\partial x_j}\rangle.
\end{equation*}
For simplicity, we denote $\frac{\partial \phi_{0,0}}{\partial
x_i}$ by $e_i$, $g_{ij}(t,\epsilon)$ by $g_{ij}$ and
$\phi_{t,\epsilon}$ by $\phi$. It is easy to see that
\begin{equation}\label{e3.1}
    \frac{\partial}{\partial t}\mid _{t=0,\epsilon=0}g_{ij}=\langle \overline{\nabla}_{e_i}\textbf{X}, e_j\rangle
    +\langle e_i,\overline{\nabla}_{e_j}\textbf{X}\rangle,
\end{equation}
\begin{equation}\label{e3.2}
    \frac{\partial}{\partial\epsilon}\mid _{t=0,\epsilon=0}g_{ij}=\langle \overline{\nabla}_{e_i}\textbf{Y}, e_j\rangle
    +\langle e_i,\overline{\nabla}_{e_j}\textbf{Y}\rangle,
\end{equation}
and
\begin{eqnarray}\label{e3.3}
\frac{\partial^2}{\partial t\partial \epsilon}\mid
_{t=0,\epsilon=0}g_{ij}
   &=& \langle\overline{\nabla}_{e_i}\textbf{Z}+K(\textbf{Y},e_i)\textbf{X}, e_j\rangle
        +\langle e_i,\overline{\nabla}_{e_j}\textbf{Z}+K(\textbf{Y},e_j)\textbf{X}\rangle \nonumber \\
   & & +\langle \overline{\nabla}_{e_i}\textbf{X}, \overline{\nabla}_{e_j}\textbf{Y}\rangle
              +\langle \overline{\nabla}_{e_i}\textbf{Y}, \overline{\nabla}_{e_j}\textbf{X}\rangle.
\end{eqnarray}
Here, $K$ is the curvature tensor on $M$. As in Section 2, if we
denote
\begin{equation*}
    \nu_{\beta}:=\nu_{\beta,t,\varepsilon}=\frac{\det^{(\beta+1)/2}(g)}{\omega^\beta(\partial\phi/\partial x_1, \partial\phi/\partial x_2)},
\end{equation*}
then
\begin{equation*}
    L_{\beta}(\phi_{t,\varepsilon})=\int_{\Sigma}\nu_{\beta} dx_1\wedge dx_2.
\end{equation*}
We have
\begin{eqnarray*}
  \frac{\partial}{\partial t}\nu_{\beta}
  &=& \frac{\beta+1}{2}\frac{\det^{(\beta+1)/2}(g)g^{ij}\frac{\partial}{\partial t}g_{ij}}
        {\omega^\beta(\partial\phi/\partial x_1,\partial\phi/\partial x_2)} \\
  & & -\beta\frac{\det^{(\beta+1)/2}(g)[\omega(\overline{\nabla}_{\partial\phi/\partial x_1}\partial\phi/\partial t, \partial\phi/\partial x_2)+\omega(\partial\phi/\partial x_1, \overline{\nabla}_{\partial\phi/\partial x_2}\partial\phi/\partial t)]}{\omega^{\beta+1}(\partial\phi/\partial x_1, \partial\phi/\partial x_2)}
\end{eqnarray*}
and
\begin{eqnarray*}
   & & \frac{\partial^2}{\partial\epsilon\partial t}\mid _{t=0,\epsilon=0} \nu_{\beta}\nonumber\\
   &=& \frac{(\beta+1)^2}{4}\frac{g^{ij}\frac{\partial}{\partial t}\mid _{t=0,\epsilon=0}g_{ij}g^{kl}\frac{\partial}{\partial \epsilon}\mid _{t=0,\epsilon=0}g_{kl}}{\cos^{\beta}\alpha}
     + \frac{\beta+1}{2}\frac{g^{ij}\frac{\partial^2}{\partial t\partial \epsilon}\mid _{t=0,\epsilon=0}g_{ij}}{\cos^{\beta}\alpha}\\
   & & -\frac{\beta+1}{2}\frac{\frac{\partial}{\partial t}\mid _{t=0,\epsilon=0}g_{ij}\frac{\partial}{\partial \epsilon}\mid _{t=0,\epsilon=0}g_{ij}}{\cos^{\beta}\alpha}\nonumber\\
   & &  -\frac{\beta(\beta+1)}{2}\frac{g^{ij}\frac{\partial}{\partial t}\mid _{t=0,\epsilon=0}g_{ij}[\omega(\overline\nabla_{e_1}\frac{\partial \phi}{\partial \epsilon}\mid_{t=0,\epsilon=0},e_2)+\omega(e_1,\overline\nabla_{e_2}\frac{\partial \phi}{\partial \epsilon}\mid_{t=0,\epsilon=0})]}{\cos^{\beta+1}\alpha}\nonumber\\
   & &  -\frac{\beta(\beta+1)}{2}\frac{g^{ij}\frac{\partial}{\partial \epsilon}\mid _{t=0,\epsilon=0}g_{ij}[\omega(\overline\nabla_{e_1}\frac{\partial \phi}{\partial t}\mid_{t=0,\epsilon=0},e_2)+\omega(e_1,\overline\nabla_{e_2}\frac{\partial \phi}{\partial t}\mid_{t=0,\epsilon=0})]}{\cos^{\beta+1}\alpha}\nonumber\\
   & & -\beta\frac{\omega(\overline\nabla_{e_1}\frac{\partial \phi}{\partial t}\mid_{t=0,\epsilon=0},\overline\nabla_{e_2}\frac{\partial \phi}{\partial \epsilon}\mid_{t=0,\epsilon=0})+\omega(\overline\nabla_{e_1}\frac{\partial \phi}{\partial \epsilon}\mid_{t=0,\epsilon=0},\overline\nabla_{e_2}\frac{\partial \phi}{\partial t}\mid_{t=0,\epsilon=0})}{\cos^{\beta+1}\alpha}\nonumber\\
   & & -\beta\frac{\omega(\overline\nabla_{e_1}\frac{\partial^2}{\partial t\partial \epsilon}\mid _{t=0,\epsilon=0}+K(\frac{\partial \phi}{\partial \epsilon}\mid_{t=0,\epsilon=0},e_1)\frac{\partial \phi}{\partial t}\mid_{t=0,\epsilon=0},e_2)}{\cos^{\beta+1}\alpha}\nonumber\\
   & & -\beta\frac{\omega(e_1,\overline\nabla_{e_2}\frac{\partial^2}{\partial t\partial \epsilon}\mid _{t=0,\epsilon=0}+K(\frac{\partial \phi}{\partial \epsilon}\mid_{t=0,\epsilon=0},e_2)\frac{\partial \phi}{\partial t}\mid_{t=0,\epsilon=0})}{\cos^{\beta+1}\alpha}\nonumber\\
   & &  +\frac{\beta(\beta+1)}{\cos^{\beta+2}\alpha}[\omega(\overline\nabla_{e_1}\frac{\partial \phi}{\partial t}\mid_{t=0,\epsilon=0},e_2)+\omega(e_1,\overline\nabla_{e_2}\frac{\partial \phi}{\partial t}\mid_{t=0,\epsilon=0})]\\
   & &  \cdot [\omega(\overline\nabla_{e_1}\frac{\partial \phi}{\partial \epsilon}\mid_{t=0,\epsilon=0},e_2)+\omega(e_1,\overline\nabla_{e_2}\frac{\partial \phi}{\partial \epsilon}\mid_{t=0,\epsilon=0})]\nonumber.
\end{eqnarray*}
By (\ref{e3.1}), (\ref{e3.2}) and (\ref{e3.3}), we have
\begin{eqnarray}\label{e3.4}
   & & \frac{\partial^2}{\partial\epsilon\partial t}\mid _{t=0,\epsilon=0} \nu_{\beta}\nonumber\\
   &=& (\beta+1)^2\frac{\langle \overline{\nabla}_{e_i}\textbf{X}, e_i\rangle\langle \overline{\nabla}_{e_j}\textbf{Y}, e_j\rangle}{\cos^{\beta}\alpha}
     + (\beta+1)\frac{\langle \overline{\nabla}_{e_i}\textbf{Z}, e_i\rangle}{\cos^{\beta}\alpha}\nonumber\\
   & & + (\beta+1)\frac{\langle \overline{\nabla}_{e_i}\textbf{X}, \overline{\nabla}_{e_i}\textbf{Y}\rangle}{\cos^{\beta}\alpha} -(\beta+1)\frac{K(\textbf{X},e_i,\textbf{Y},e_i)}{\cos^{\beta}\alpha}\nonumber\\
   & &  -\frac{\beta+1}{2}\frac{(\langle \overline{\nabla}_{e_i}\textbf{X}, e_j\rangle
    +\langle e_i,\overline{\nabla}_{e_j}\textbf{X}\rangle)(\langle \overline{\nabla}_{e_i}\textbf{Y}, e_j\rangle
    +\langle e_i,\overline{\nabla}_{e_j}\textbf{Y}\rangle)}{\cos^{\beta}\alpha}\nonumber\\
   & &  -\beta(\beta+1)\frac{\langle \overline{\nabla}_{e_i}\textbf{X}, e_i\rangle[\omega(\overline\nabla_{e_1}\textbf{Y},e_2)+\omega(e_1,\overline\nabla_{e_2}\textbf{Y})]}{\cos^{\beta+1}\alpha}\nonumber\\
   & &  -\beta(\beta+1)\frac{\langle \overline{\nabla}_{e_i}\textbf{Y}, e_i\rangle[\omega(\overline\nabla_{e_1}\textbf{X},e_2)+\omega(e_1,\overline\nabla_{e_2}\textbf{X})]}{\cos^{\beta+1}\alpha}\nonumber\\
   & & -\beta\frac{\omega(\overline\nabla_{e_1}\textbf{X},\overline\nabla_{e_2}\textbf{Y})
   +\omega(\overline\nabla_{e_1}\textbf{Y},\overline\nabla_{e_2}\textbf{X})}{\cos^{\beta+1}\alpha}\nonumber\\
   & & -\beta\frac{\omega(\overline\nabla_{e_1}\textbf{Z},e_2)+\omega(e_1,\overline\nabla_{e_2}\textbf{Z})}{\cos^{\beta+1}\alpha}\nonumber\\
   & & -\beta\frac{\omega(K(\textbf{Y},e_1)\textbf{X},e_2)+\omega(e_1,K(\textbf{Y},e_2)\textbf{X})}{\cos^{\beta+1}\alpha}\nonumber\\
   & &  +\beta(\beta+1)\frac{[\omega(\overline\nabla_{e_1}\textbf{X},e_2)+\omega(e_1,\overline\nabla_{e_2}\textbf{X})]
         [\omega(\overline\nabla_{e_1}\textbf{Y},e_2)+\omega(e_1,\overline\nabla_{e_2}\textbf{Y})]}{\cos^{\beta+2}\alpha}.
\end{eqnarray}
Since $\Sigma$ is a $\beta$-symplectic critical surface, by the first
variation formula (\ref{e2.4}), we see that
\begin{equation}\label{e3.5}
     (\beta+1)\int_{\Sigma}\frac{\langle \overline{\nabla}_{e_i}\textbf{Z}, e_i\rangle}{\cos^\beta\alpha}d\mu
-\beta\int_{\Sigma}\frac{\omega(\overline\nabla_{e_1}
\textbf{Z}, e_2)+ \omega(e_1, \overline\nabla_{e_2}
\textbf{Z})}{\cos^{\beta+1}\alpha}d\mu=0.
\end{equation}
Notice that,
\begin{eqnarray}\label{e3.6}
  & & 2\langle \overline{\nabla}_{e_i}\textbf{X}, \overline{\nabla}_{e_i}\textbf{Y}\rangle
   -(\langle \overline{\nabla}_{e_i}\textbf{X}, e_j\rangle
    +\langle e_i,\overline{\nabla}_{e_j}\textbf{X}\rangle)(\langle \overline{\nabla}_{e_i}\textbf{Y}, e_j\rangle
    +\langle e_i,\overline{\nabla}_{e_j}\textbf{Y}\rangle)\nonumber \\
  &=&  2\langle \overline{\nabla}^{\perp}_{e_i}\textbf{X}, \overline{\nabla}^{\perp}_{e_i}\textbf{Y}\rangle
      -(\langle e_i,\overline{\nabla}_{e_j}\textbf{X}\rangle\langle e_j,\overline{\nabla}_{e_i}\textbf{Y}\rangle
       +\langle e_i,\overline{\nabla}_{e_j}\textbf{Y}\rangle\langle e_j,\overline{\nabla}_{e_i}\textbf{X}\rangle),
\end{eqnarray}
and
\begin{equation}\label{e3.7}
    \langle \overline{\nabla}_{e_i}\textbf{X}, e_i\rangle=div_{\Sigma}\textbf{X}.
\end{equation}
Putting (\ref{e3.5}), (\ref{e3.6}) and (\ref{e3.7}) into
(\ref{e3.4}), we obtain the identity in the lemma. \hfill Q.E.D.

\vspace{.1in}

Assume now that $\textbf{X}=\textbf{Y}$ is a normal vector field,
then we have
\begin{equation}\label{e3.8}
    div_{\Sigma}\textbf{X}=-\langle\textbf{X},\textbf{H}\rangle,
\end{equation}
\begin{equation}\label{e3.9}
    \langle e_i,\overline{\nabla}_{e_j}\textbf{X}\rangle\langle e_j,\overline{\nabla}_{e_i}\textbf{X}\rangle
    =\langle\textbf{X},\textbf{A}(e_i,e_j)\rangle^2,
\end{equation}
and
\begin{eqnarray}\label{e3.10}
   & & -\beta\int_{\Sigma}\frac{\omega(K(\textbf{X},e_1)\textbf{X},e_2)+\omega(e_1,K(\textbf{X},e_2)\textbf{X})}{\cos^{\beta+1}\alpha}d\mu
       -2\beta\int_{\Sigma}\frac{\omega(\overline\nabla_{e_1}\textbf{X},\overline\nabla_{e_2}\textbf{X})}{\cos^{\beta+1}\alpha}d\mu\nonumber\\
   &=& \beta\int_{\Sigma}\frac{\langle Je_2,K(\textbf{X},e_1)\textbf{X})-\langle Je_1,K(\textbf{X},e_2)\textbf{X})}
                   {\cos^{\beta+1}\alpha}d\mu\nonumber\\
   & & -\beta\int_{\Sigma}\frac{e_1(\omega(\textbf{X},\overline\nabla_{e_2}\textbf{X}))-\omega(\textbf{X},
\overline\nabla_{e_1}\overline\nabla_{e_2} \textbf{X})}{\cos^{\beta+1}\alpha}d\mu\nonumber\\
   & & -\beta\int_{\Sigma}\frac{e_2(\omega(\overline\nabla_{e_1}\textbf{X},\textbf{X}))
         -\omega(\overline\nabla_{e_2}\overline\nabla_{e_1}\textbf{X},\textbf{X})}{\cos^{\beta+1}\alpha}d\mu\nonumber\\
   &=& \beta\int_{\Sigma}\frac{K(\textbf{X},e_1,Je_2,\textbf{X})+K(\textbf{X},e_2,\textbf{X},Je_1)+K(e_1,e_2,J\textbf{X},\textbf{X})}
                   {\cos^{\beta+1}\alpha}d\mu\nonumber\\
   & & -\beta(\beta+1)\int_{\Sigma}\frac{\omega(\textbf{X},\overline\nabla_{e_2}\textbf{X})\nabla_{e_1}\cos\alpha
              +\omega(\overline\nabla_{e_1}\textbf{X},\textbf{X}))\nabla_{e_2}\cos\alpha}{\cos^{\beta+2}\alpha}d\mu.
\end{eqnarray}
In the last step, we used integration by parts. By the first
Bianchi identity, we have
\begin{eqnarray}\label{e3.11}
  & & K(\textbf{X},e_1,Je_2,\textbf{X})+K(\textbf{X},e_2,\textbf{X},Je_1)+K(e_1,e_2,J\textbf{X}),\textbf{X}\nonumber \\
  &=& K(\textbf{X},e_1,J\textbf{X},e_2)+K(\textbf{X},e_2,e_1,J\textbf{X})+K(\textbf{X},J\textbf{X},e_2,e_1)=0.
\end{eqnarray}
Therefore, by (\ref{e3.8}), (\ref{e3.9}), (\ref{e3.10}) and
(\ref{e3.11}), we obtain:

\begin{corollary}\label{cor4.2}
If $\textbf{X}$ is a normal variation, then
\begin{eqnarray}\label{2nd2}
   & & \frac{\partial^2}{\partial t^2}\mid _{t=0} L_{\beta}(\phi_{t}):= II_{\beta}(\textbf{X})\nonumber\\
   &=& (\beta+1)\int_{\Sigma}\frac{|\overline{\nabla}_{e_i}^{\perp}\textbf{X}|^2}{\cos^{\beta}\alpha} d\mu
       -(\beta+1)\int_{\Sigma}\frac{\langle\textbf{X},\textbf{A}(e_i,e_j)\rangle^2}{\cos^{\beta}\alpha} d\mu
         -(\beta+1)\int_{\Sigma}\frac{K(\textbf{X},e_i,\textbf{X},e_i)}{\cos^{\beta}\alpha}d\mu\nonumber\\
   & & +(\beta+1)^2\int_{\Sigma}\frac{\langle \textbf{X},\textbf{H}\rangle^2}{\cos^{\beta}\alpha} d\mu
        +2\beta(\beta+1)\int_{\Sigma}\frac{\langle \textbf{X},\textbf{H}\rangle
             (\omega(\overline{\nabla}_{e_1}\textbf{X}, e_2)+\omega(e_1,\overline{\nabla}_{e_2}\textbf{X}))}{\cos^{\beta+1}\alpha}d\mu\nonumber\\
   & &  -\beta(\beta+1)\int_{\Sigma}\frac{\omega(\textbf{X},\overline{\nabla}_{e_2}\textbf{X})\nabla_{e_1}\cos\alpha
              +\omega(\overline{\nabla}_{e_1}\textbf{X},\textbf{X})\nabla_{e_2}\cos\alpha}{\cos^{\beta+2}\alpha} d\mu\nonumber\\
   & & +\beta(\beta+1)\int_{\Sigma}\frac{(\omega(\overline{\nabla}_{e_1}\textbf{X}, e_2)
          +\omega(e_1,\overline{\nabla}_{e_2}\textbf{X}))^2}{\cos^{\beta+2}\alpha}d\mu.
\end{eqnarray}
\end{corollary}

\vspace{.1in}

Next, we will follow Micallef-Wolfson's idea (\cite{MW}) to give
another version of the second variation formula so that we can
have some rigidity results as its corollaries. For the background
material of the geometry of a surface in a four-manifold, we refer
to Section 3 of Micallef-Wolfson's paper (\cite{MW}).

\vspace{.1in}

Set $\textbf{X}=x_3e_3+x_4e_4$, then
\begin{eqnarray*}
  \omega(\textbf{X},\overline{\nabla}_{e_2}\textbf{X})
   &=& \langle J\textbf{X},\overline{\nabla}_{e_2}\textbf{X}\rangle \\
   &=& \langle x_3Je_3+x_4Je_4,x_{32}e_3+x_{42}e_4-x_3h^3_{2i}e_i-x_4h^4_{2l}e_l\rangle \\
   &=& x_3x_{42}\langle Je_3,e_4\rangle+x_4x_{32}\langle Je_4,e_3\rangle \\
   & & +x_3(x_3h^3_{12}+x_4h^4_{12})\langle e_3,Je_1\rangle+x_3(x_3h^3_{22}+x_4h^4_{22})\langle e_3,Je_2\rangle  \\
   & & +x_4(x_3h^3_{12}+x_4h^4_{12})\langle e_4,Je_1\rangle+x_4(x_3h^3_{22}+x_4h^4_{22})\langle e_4,Je_2\rangle.
\end{eqnarray*}
Notice that (see (\ref{e15}))
\begin{equation}\label{e3.14}
\langle Je_1,e_2\rangle=\langle Je_3,e_4\rangle=\cos\alpha,
\end{equation}
\begin{equation}\label{e3.15}
\langle Je_1,e_3\rangle=-\langle Je_2,e_4\rangle,
\end{equation}
\begin{equation}\label{e3.16}
\langle Je_1,e_4\rangle=\langle Je_2,e_3\rangle.
\end{equation}
Thus, we have
\begin{eqnarray}\label{e3.17}
  \omega(\textbf{X},\overline{\nabla}_{e_2}\textbf{X})
  &=& (x_3x_{42}-x_4x_{32})\cos\alpha\nonumber \\
  & &  +(x_3^2h^3_{12}+x_3x_4h^4_{12}-x_3x_4h^3_{22}-x_4^2h^4_{22})\langle e_3,Je_1\rangle\nonumber \\
  & &  +(x_3^2h^3_{22}+x_3x_4h^4_{22}+x_3x_4h^3_{12}+x_4^2h^4_{12})\langle e_3,Je_2\rangle,
\end{eqnarray}
and
\begin{eqnarray}\label{e3.18}
  \omega(\textbf{X},\overline{\nabla}_{e_1}\textbf{X})
  &=& (x_3x_{41}-x_4x_{31})\cos\alpha\nonumber \\
  & &  +(x_3^2h^3_{11}+x_3x_4h^4_{11}-x_3x_4h^3_{12}-x_4^2h^4_{12})\langle e_3,Je_1\rangle\nonumber \\
  & &  +(x_3^2h^3_{12}+x_3x_4h^4_{12}+x_3x_4h^3_{11}+x_4^2h^4_{11})\langle e_3,Je_2\rangle,
\end{eqnarray}
Similarly,
\begin{eqnarray}\label{e3.19}
  \omega(\overline{\nabla}_{e_1}\textbf{X}, e_2)
          +\omega(e_1,\overline{\nabla}_{e_2}\textbf{X})
  &=& -\cos\alpha(x_3H^3+x_4H^4)\\
  & &  +(x_{32}+x_{41})\langle Je_1, e_3\rangle+(x_{42}-x_{31})\langle Je_2, e_3\rangle.\nonumber
\end{eqnarray}
By (\ref{e3.19}), we get that
\begin{eqnarray}\label{e3.20}
   & & 2\beta(\beta+1)\int_{\Sigma}\frac{\langle \textbf{X},\textbf{H}\rangle
             (\omega(\overline{\nabla}_{e_1}\textbf{X}, e_2)+\omega(e_1,\overline{\nabla}_{e_2}\textbf{X}))}{\cos^{\beta+1}\alpha}d\mu\nonumber\\
   & & +\beta(\beta+1)\int_{\Sigma}\frac{(\omega(\overline{\nabla}_{e_1}\textbf{X}, e_2)
                 +\omega(e_1,\overline{\nabla}_{e_2}\textbf{X}))^2}{\cos^{\beta+2}\alpha}d\mu\nonumber\\
   &=& -\beta(\beta+1)\int_{\Sigma}\frac{(x_3H^3+x_4H^4)^2}{\cos^{\beta}\alpha}d\mu
     +\beta(\beta+1)\int_{\Sigma}\frac{(x_{32}+x_{41})^2\langle Je_1, e_3\rangle^2}{\cos^{\beta+2}\alpha}d\mu\nonumber\\
   & & +\beta(\beta+1)\int_{\Sigma}\frac{(x_{42}-x_{31})^2\langle Je_2, e_3\rangle^2}{\cos^{\beta+2}\alpha}d\mu\nonumber\\
   & & +2\beta(\beta+1)\int_{\Sigma}\frac{(x_{32}+x_{41})(x_{42}-x_{31})\langle Je_1, e_3\rangle\langle Je_2, e_3\rangle}{\cos^{\beta+2}\alpha}d\mu.
\end{eqnarray}
Now $II_{\beta}(\textbf{X})$ can be written as
\begin{eqnarray}\label{e3.21}
  II_{\beta}(\textbf{X})
   &=& (\beta+1)\int_{\Sigma}\frac{x_{3i}^2+x_{4i}^2}{\cos^{\beta}\alpha} d\mu
         -(\beta+1)\int_{\Sigma}\frac{x_3^2K_{3i3i}+x_4^2K_{4i4i}+2x_3x_4K_{3i4i}}{\cos^{\beta}\alpha}d\mu\nonumber\\
   & & -2(\beta+1)\int_{\Sigma}\frac{(x_3h^3_{12}+x_4h^4_{12})^2-(x_3h^3_{11}+x_4h^4_{11})(x_3h^3_{22}+x_4h^4_{22})}{\cos^{\beta}\alpha}d\mu\nonumber\\
   & &  +\beta(\beta+1)\int_{\Sigma}\frac{(x_{32}+x_{41})^2\langle Je_1, e_3\rangle^2+(x_{42}-x_{31})^2\langle Je_2, e_3\rangle^2}{\cos^{\beta+2}\alpha}d\mu\nonumber\\
   & & +2\beta(\beta+1)\int_{\Sigma}\frac{(x_{32}+x_{41})(x_{42}-x_{31})\langle Je_1, e_3\rangle\langle Je_2, e_3\rangle}{\cos^{\beta+2}\alpha}d\mu\nonumber\\
   & & +III_{\beta}(\textbf{X}),
\end{eqnarray}
where
\begin{equation*}
    III_{\beta}(\textbf{X})=-\beta(\beta+1)\int_{\Sigma}\frac{\omega(\textbf{X},\overline{\nabla}_{e_2}\textbf{X})\nabla_{e_1}\cos\alpha
              +\omega(\overline{\nabla}_{e_1}\textbf{X},\textbf{X})\nabla_{e_2}\cos\alpha}{\cos^{\beta+2}\alpha} d\mu
\end{equation*}

In addition, if we set
$\textbf{Y}=-J_{\nu}\textbf{X}=x_4e_3-x_3e_4$, then by
(\ref{e3.21}), we get that
\begin{eqnarray}\label{e3.22}
  & & II_{\beta}(\textbf{X})+II_{\beta}(\textbf{Y})\nonumber\\
   &=& 2(\beta+1)\int_{\Sigma}\frac{x_{3i}^2+x_{4i}^2}{\cos^{\beta}\alpha} d\mu
         -(\beta+1)\int_{\Sigma}\frac{|\textbf{X}|^2(K_{3i3i}+K_{4i4i})}{\cos^{\beta}\alpha}d\mu\nonumber\\
   & & -2(\beta+1)\int_{\Sigma}\frac{|\textbf{X}|^2[(h^3_{12})^2+(h^4_{12})^2-h^3_{11}h^3_{22}-h^4_{11}h^4_{22}]}{\cos^{\beta}\alpha}d\mu\nonumber\\
   & &  +\beta(\beta+1)\int_{\Sigma}\frac{[(x_{32}+x_{41})^2+(x_{42}-x_{31})^2][\langle Je_1, e_3\rangle^2+\langle Je_2, e_3\rangle^2]}{\cos^{\beta+2}\alpha}d\mu\nonumber\\
   & & +III_{\beta}(\textbf{X})+III_{\beta}(\textbf{Y}).
\end{eqnarray}
Using (\ref{e3.17}) and (\ref{e3.18}), it is easy to check that
\begin{eqnarray}\label{e3.23}
  III_{\beta}(\textbf{X})+III_{\beta}(\textbf{Y})
   &=& 2(\beta+1)\int_{\Sigma}(x_3x_{42}-x_4x_{32})\nabla_{e_1}\frac{1}{\cos^{\beta}\alpha} d\mu\\
   & & -2(\beta+1)\int_{\Sigma}(x_3x_{41}-x_4x_{31})\nabla_{e_2}\frac{1}{\cos^{\beta}\alpha} d\mu\nonumber\\
   & & -\beta(\beta+1)\int_{\Sigma}\frac{|\textbf{X}|^2(h^3_{12}-h^4_{22})\langle e_3,Je_1\rangle\nabla_{e_1}\cos\alpha}{\cos^{\beta+2}\alpha}d\mu\nonumber\\
   & & -\beta(\beta+1)\int_{\Sigma}\frac{|\textbf{X}|^2(h^3_{22}+h^4_{12})\langle e_4,Je_1\rangle\nabla_{e_1}\cos\alpha}{\cos^{\beta+2}\alpha}d\mu\nonumber\\
   & & +\beta(\beta+1)\int_{\Sigma}\frac{|\textbf{X}|^2(h^3_{11}-h^4_{12})\langle e_3,Je_1\rangle\nabla_{e_2}\cos\alpha}{\cos^{\beta+2}\alpha}d\mu\nonumber\\
   & & +\beta(\beta+1)\int_{\Sigma}\frac{|\textbf{X}|^2(h^3_{12}+h^4_{11})\langle e_4,Je_1\rangle\nabla_{e_2}\cos\alpha}{\cos^{\beta+2}\alpha}d\mu.\nonumber
\end{eqnarray}
Integrating the first two terms in (\ref{e3.23}) by parts, we get
that
\begin{eqnarray*}
   & & 2(\beta+1)\int_{\Sigma}(x_3x_{42}-x_4x_{32})\nabla_{e_1}\frac{1}{\cos^{\beta}\alpha} d\mu -2(\beta+1)\int_{\Sigma}(x_3x_{41}-x_4x_{31})\nabla_{e_2}\frac{1}{\cos^{\beta}\alpha} d\mu\nonumber\\
   &=& 2(\beta+1)\int_{\Sigma}\frac{2(x_{32}x_{41}-x_{31}x_{42})+x_3(x_{412}-x_{421})+x_4(x_{321}-x_{312})}{\cos^{\beta}\alpha} d\mu\nonumber\\
   &=& 4(\beta+1)\int_{\Sigma}\frac{x_{32}x_{41}-x_{31}x_{42}}{\cos^{\beta}\alpha} d\mu+2(\beta+1)\int_{\Sigma}\frac{|\textbf{X}|^2R_{3412}}{\cos^{\beta}\alpha} d\mu,
\end{eqnarray*}
where $R_{3412}$ is the normal curvature of $\Sigma$ in $M$. Here,
we have used the fact that (see (3.9) of (\cite{MW}))
\begin{equation*}
-x_3R_{3412}=x_{421}-x_{412}, \ \ \ x_4R_{3412}=x_{321}-x_{312}.
\end{equation*}
We will compute the last four terms of (\ref{e3.23}) at each
point, so we can choose an orthonormal frame $\{e_1,e_2,e_3,e_4\}$
around the point so that $J$ takes the form
\begin{eqnarray}\label{e3.24}
J=\left (\begin{array}{clcr} 0 &\cos\alpha &\sin\alpha &0 \\
-\cos\alpha &0 &0 &-\sin\alpha\\
-\sin\alpha &0 &0 &\cos\alpha\\
0 &\sin\alpha &-\cos\alpha &0 \end{array}\right).
\end{eqnarray}
By the equation of a $\beta$-symplectic critical surface
(\ref{equation2}), we have
\begin{equation*}
    H^4=\beta\frac{\sin^2\alpha}{\cos^2\alpha}\nabla_{e_1}\alpha, \ \ H^3=\beta\frac{\sin^2\alpha}{\cos^2\alpha}\nabla_{e_2}\alpha.
\end{equation*}
Using (\ref{e.alpha}), we have
\begin{eqnarray*}
   & & -(h^3_{12}-h^4_{22})\nabla_{e_1}\cos\alpha+(h^3_{11}-h^4_{12})\nabla_{e_2}\cos\alpha\\
   &=& \sin\alpha(h^3_{12}+h^4_{11}-H^4)\nabla_{e_1}\alpha+\sin\alpha(h^3_{22}+h^4_{12}-H^3)\nabla_{e_2}\alpha\nonumber\\
   &= & -\sin\alpha(\nabla_{e_1}\alpha+\beta\frac{\sin^2\alpha}{\cos^2\alpha}\nabla_{e_1}\alpha)\nabla_{e_1}\alpha
    -\sin\alpha(\nabla_{e_2}\alpha+\beta\frac{\sin^2\alpha}{\cos^2\alpha}\nabla_{e_2}\alpha)\nabla_{e_2}\alpha\nonumber\\
   &=& -\frac{\sin\alpha(\cos^2\alpha+\beta\sin^2\alpha)}{\cos^2\alpha}|\nabla\alpha|^2,
\end{eqnarray*}
and
\begin{equation*}
 -(h^3_{22}+h^4_{12})\nabla_{e_1}\cos\alpha+(h^3_{12}+h^4_{11})\nabla_{e_2}\cos\alpha
   = -\sin\alpha\nabla_{e_2}\alpha\nabla_{e_1}\alpha+\sin\alpha\nabla_{e_1}\alpha\nabla_{e_2}\alpha=0.
\end{equation*}
By the Ricci equation,
\begin{equation*}
    R_{3412}=K_{3412}+h^{3}_{1k}h^4_{2k}-h^3_{2k}h^4_{1k},
\end{equation*}
we obtain that
\begin{eqnarray}\label{e3.26}
  III_{\beta}(\textbf{X})+III_{\beta}(\textbf{Y})
   &=& 4(\beta+1)\int_{\Sigma}\frac{x_{32}x_{41}-x_{31}x_{42}}{\cos^{\beta}\alpha} d\mu+2(\beta+1)\int_{\Sigma}\frac{|\textbf{X}|^2K_{3412}}{\cos^{\beta}\alpha} d\mu\nonumber\\
   & & +2(\beta+1)\int_{\Sigma}\frac{|\textbf{X}|^2(h^{3}_{1k}h^4_{2k}-h^3_{2k}h^4_{1k})}{\cos^{\beta}\alpha} d\mu\nonumber\\
   & & -\beta(\beta+1)\int_{\Sigma}\frac{|\textbf{X}|^2\sin^2\alpha(\cos^2\alpha+\beta\sin^2\alpha)}{\cos^{\beta+4}\alpha}|\nabla\alpha|^2d\mu.
\end{eqnarray}
Putting (\ref{e3.26}) into (\ref{e3.22}) yields
\begin{eqnarray*}
   & & II_{\beta}(\textbf{X})+II_{\beta}(\textbf{Y})\\
   &=& 2(\beta+1)\int_{\Sigma}\frac{(x_{32}+x_{41})^2+(x_{42}-x_{31})^2}{\cos^{\beta}\alpha} d\mu\\
   & & -(\beta+1)\int_{\Sigma}\frac{|\textbf{X}|^2(K_{3i3i}+K_{4i4i}-2K_{1234})}{\cos^{\beta}\alpha}d\mu\nonumber\\
   & &  +\beta(\beta+1)\int_{\Sigma}\frac{[(x_{32}+x_{41})^2+(x_{42}-x_{31})^2]\sin^2\alpha}{\cos^{\beta+2}\alpha}d\mu\nonumber\\
   & & -\beta(\beta+1)\int_{\Sigma}\frac{|\textbf{X}|^2\sin^2\alpha(\cos^2\alpha+\beta\sin^2\alpha)}{\cos^{\beta+4}\alpha}|\nabla\alpha|^2d\mu\nonumber\\   & & -2(\beta+1)\int_{\Sigma}\frac{|\textbf{X}|^2[(h^3_{12})^2+(h^4_{12})^2-h^3_{11}h^3_{22}-h^4_{11}h^4_{22}-h^{3}_{1k}h^4_{2k}+h^3_{2k}h^4_{1k}]}{\cos^{\beta}\alpha}d\mu.
\end{eqnarray*}
By Lemma 3.2 of \cite{HL}, it is easy to check that
\begin{equation*}
    K_{3i3i}+K_{4i4i}-2K_{1234}=2K\sin^2\alpha,
\end{equation*}
where $R$ is the scalar curvature of $M$. By (3.16) of \cite{MW},
we have
\begin{equation*}
    (x_{32}+x_{41})^2+(x_{42}-x_{31})^2=|\bar\partial \textbf{X}|^2.
\end{equation*}
Using the fact that $\Sigma$ is a $\beta$-symplectic critical surface
again, we can see that
\begin{eqnarray*}
   & & (h^3_{12})^2+(h^4_{12})^2-h^3_{11}h^3_{22}-h^4_{11}h^4_{22}-h^{3}_{1k}h^4_{2k}+h^3_{2k}h^4_{1k}\\
   &=& (h^4_{11}+h^3_{12})(h^3_{12}-h^4_{22})+(h^4_{12}+h^3_{22})(h^4_{12}-h^3_{11})\\
   &=& -\nabla_{e_1}\alpha(-\nabla_{e_1}\alpha-H^4)-\nabla_{e_2}\alpha(-\nabla_{e_2}\alpha-H^3)\nonumber\\
   &=&  \nabla_{e_1}\alpha(\nabla_{e_1}\alpha+\beta\frac{\sin^2\alpha}{\cos^2\alpha}\nabla_{e_1}\alpha)+\nabla_{e_2}\alpha(\nabla_{e_2}\alpha+\beta\frac{\sin^2\alpha}{\cos^2\alpha}\nabla_{e_2}\alpha)\nonumber\\   &=& \frac{\cos^2\alpha+\beta\sin^2\alpha}{\cos^2\alpha}|\nabla\alpha|^2
\end{eqnarray*}
Therefore, we can conclude that
\begin{eqnarray*}
  II_{\beta}(\textbf{X})+II_{\beta}(\textbf{Y})
   &=& -2(\beta+1)\int_{\Sigma}\frac{K|\textbf{X}|^2\sin^2\alpha}{\cos^{\beta}\alpha}d\mu +2(\beta+1)\int_{\Sigma}\frac{|\bar\partial \textbf{X}|^2}{\cos^{\beta}\alpha} d\mu\\
   & &  +\beta(\beta+1)\int_{\Sigma}\frac{|\bar\partial \textbf{X}|^2\sin^2\alpha}{\cos^{\beta+2}\alpha}d\mu\nonumber\\
   & & -\beta(\beta+1)\int_{\Sigma}\frac{|\textbf{X}|^2\sin^2\alpha(\cos^2\alpha+\beta\sin^2\alpha)}{\cos^{\beta+4}\alpha}|\nabla\alpha|^2d\mu\nonumber\\   & & -2(\beta+1)\int_{\Sigma}\frac{|\textbf{X}|^2(\cos^2\alpha+\beta\sin^2\alpha)}{\cos^{\beta+2}\alpha}|\nabla\alpha|^2d\mu\\
   &=& -2(\beta+1)\int_{\Sigma}\frac{K|\textbf{X}|^2\sin^2\alpha}{\cos^{\beta}\alpha}d\mu +(\beta+1)\int_{\Sigma}\frac{|\bar\partial \textbf{X}|^2(2\cos^2\alpha+\beta\sin^2\alpha)}{\cos^{\beta+2}\alpha} d\mu\\
   & & -(\beta+1)\int_{\Sigma}\frac{(2\cos^2\alpha+\beta\sin^2\alpha)(\cos^2\alpha+\beta\sin^2\alpha)}{\cos^{\beta+4}\alpha}|\textbf{X}|^2|\nabla\alpha|^2d\mu.
\end{eqnarray*}
This gives the following theorem:

\begin{theorem}\label{thm4.3}
If we choose $\textbf{X}=x_3e_3+x_4e_4$ and
$\textbf{Y}=-J_{\nu}\textbf{X}=x_4e_3-x_3e_4$, then the second
variation of the functional $L_{\beta}$ of a $\beta$-symplectic
critical surface is
\begin{eqnarray}\label{2nd3}
   & & II_{\beta}(\textbf{X})+II_{\beta}(\textbf{Y})\nonumber\\
   &=& -2(\beta+1)\int_{\Sigma}\frac{K|\textbf{X}|^2\sin^2\alpha}{\cos^{\beta}\alpha}d\mu +(\beta+1)\int_{\Sigma}\frac{|\bar\partial \textbf{X}|^2(2\cos^2\alpha+\beta\sin^2\alpha)}{\cos^{\beta+2}\alpha} d\mu\nonumber\\
   & & -(\beta+1)\int_{\Sigma}\frac{(2\cos^2\alpha+\beta\sin^2\alpha)(\cos^2\alpha+\beta\sin^2\alpha)}{\cos^{\beta+4}\alpha}|\textbf{X}|^2|\nabla\alpha|^2d\mu.
\end{eqnarray}
\end{theorem}

As applications of the stability inequality (\ref{2nd3}), we can
obtain some rigidity results for stable $\beta$-symplectic
critical surfaces.

\vspace{.1in}

\begin{corollary}\label{cor4.4}
Let $M$ be a K\"ahler surface with positive scalar curvature $K$.
If $\Sigma$ is a stable $\beta$-symplectic critical surface in $M$
with $\beta\geq0$, whose normal bundle admits a nontrivial section
$\textbf{X}$ with
\begin{equation*}
    \frac{|\bar\partial \textbf{X}|^2}{|\textbf{X}|^2}\leq\frac{\cos^2\alpha+\beta\sin^2\alpha}{\cos^2\alpha}|\nabla\alpha|^2,
\end{equation*}
then $\Sigma$ is a holomorphic curve.
\end{corollary}

\textbf{Proof:} If $\Sigma$ is a stable $\beta$-symplectic
critical surface in $M$, we have
\begin{equation*}
    II_{\beta}(\textbf{X})+II_{\beta}(\textbf{Y})\geq0,
\end{equation*}
where $\textbf{Y}=-J_{\nu}\textbf{X}$. By Theorem \ref{thm4.3}
with $\beta\geq0$, we obtain
\begin{equation*}
    2(\beta+1)\int_{\Sigma}\frac{K|\textbf{X}|^2\sin^2\alpha}{\cos^{\beta}\alpha}d\mu\leq0,
\end{equation*}
Since $\beta\geq0$ and $K>0$, by the last inequality, we must have
$\sin\alpha\equiv0$, that is, $\Sigma$ is a holomorphic curve.
\hfill Q.E.D.

\vspace{.1in}

\begin{corollary}\label{cor4.5}
Let $M$ be a K\"ahler surface with positive scalar curvature $K$.
If $\Sigma$ is a stable $\beta$-symplectic critical surface in $M$
with $\beta\geq0$ and $\chi(\nu)\geq g$, where $\chi(\nu)$ is the
Euler characteristic of the normal bundle $\nu$ of $\Sigma$ in $M$
and $g$ is the genus of $\Sigma$, then $\Sigma$ is a holomorphic
curve.
\end{corollary}

\textbf{Proof:} Let $L$ be the canonical line bundle over
$\Sigma$, let $H^0(\Sigma,{\mathcal O}(\nu))$ denote the space of
holomorphic sections of $\nu$, and let $H^0(\Sigma,{\mathcal
O}(\nu\otimes L))$ denote the space of holomorphic sections of
$\nu\otimes L$, then by Riemann-Roch theorem, we have
\begin{equation*}
    \dim H^0(\Sigma,{\mathcal O}(\nu))=\chi(\nu)-g+1+\dim H^0(\Sigma,{\mathcal O}(\nu\otimes L))\geq 1.
\end{equation*}
We therefore have a nontrivial holomorphic section on $\nu$ and
the corollary follows from Corollary \ref{cor4.4}. \hfill Q.E.D.

\vspace{.1in}

\section{Continuity Method and Openness}

We would like to study the family of $\beta$-symplectic critical surfaces in $[0,\infty)$. We set
\begin{equation*}
    S:=\{\beta\in[0,\infty)\mid \exists\ a \ strictly \ stable\ \beta-symplecitc\ critical\ surface\}.
\end{equation*}

\begin{theorem}\label{thm5.1}
$S$ is open on $[0,\infty)$. In other words, if
 there is a strictly stable $\beta_0$-symplectic critical surface $\Sigma_{\beta_0}$, then there is a strictly stable $\beta$-symplectic critical surface in the neighborhood of $\Sigma_{\beta_0}$ for $\beta$ sufficiently close to $\beta_0$.
\end{theorem}

Before we prove our theorem, we state  theorem $1.1$ of White \cite{White} here for the reader's convenience.
\begin{theorem}\cite{White}
Let $M$ be a smooth compact $m$-dimensional Riemannian manifold, $V$ be a smooth Euclidean vector bundle over $M$, and $\nabla$ be a smooth Riemannian connection on $V$. Let $G$ be the Banach space of $C^q$ functions $f$ that assign to each $x\in M$, $v\in V_x$, and linear map $L: Tan_x M\to V_x$ a real number $f(x, v, L)$ in such a way that $D_3f$ is also $C^q$. Let $\Gamma$ be an open subset of a Banach space and $A$ be a smooth map from $\Gamma$ to $G$. Then:\\
(1) If $u_t\in C^{2, \alpha}(M, V)$ is a differentiable one-parameter family of sections, then
$$(\frac{d}{dt})_{t=0}\int_{M} A_\gamma(x, u_t(x), \nabla u_t(x))dx=\int_{M} H(\gamma, u_0)(x)\cdot\frac{d}{dt}u_0(x)dx$$ where
$$H(\gamma, u)=-div D_3A_\gamma(x, u, \nabla u)+D_2 A_\gamma(x, u, \nabla u).$$
Furthermore, if $2\leq j\leq q$ and $u\in C^{j, \alpha}(M, V)$, then $H(\gamma, u)\in C^{j-2,\alpha}(M, V)$ and the map
$$H: \Gamma\times C^{j, \alpha}(M, V)\to C^{j-2, \alpha}(M, V)$$ is $C^{q-j}$.\\
(2) The linearization
$$J=D_2 H(\gamma_0, u_0): C^{j, \alpha}(M, V)\to C^{j-2,\alpha}(M, V)$$
is a self-adjoint second order linear partial differential operator.\\
If in addition $J$ is an elliptic operator, then\\
(3) $J$ is a Fredholm map with Fredholm index $0$.\\
(4) If $H(\gamma_0, u_0)=0$, then $u_0$ must be $C^{q, \beta}$ for every $\beta<1$.\\
(5) There is an orthonormal basis for $L^2(M, V)$ consisting of $C^q$ eigenfunctions for $J$. For every $\lambda$, only finitely many eigenfunctions have eigenvalues less than $\lambda$.
\end{theorem}

\vspace{.1in}

\textbf{Proof of Theorem \ref{thm5.1}:} It suffices to show that if $\beta_0\in S$, i.e.,
there exists a strictly stable $\beta_0$-symplecitc critical
surface $\Sigma_{\beta_0}$ in $M$, then there at least one
solution to the equation (\ref{betaequ}) which is also strictly stable
when $\beta$ is sufficiently close to $\beta_0$. We will use the
above theorem and Implicit Function Theorem.

In order to use the above theorem of White. We need equate the map $F$ with the section $u$ of normal bundle and we consider the equivalent class $[u]$ in place of $u$ ([u] denote the set of all maps $u\circ\phi$ where $\phi:\Sigma\to\Sigma$ is a diffeomorphism). See theorem $2.1$ in \cite{White}.

We have seen in Section 3 that the linearization $J_{\beta_0}=DP(\beta, [F])$ is an elliptic
operator. Therefore, by $(3)$ of the above theorem,  $J_{\beta_0}$ is a Fredholm operator
with Fredholm index 0. Since we assume that $\Sigma_{\beta_0}$ is
strictly stable so that $J_{\beta_0}$ has no normal Jacobi fields,
therefore $J_{\beta_0}$ is an isomorphism. Thus, by Implicit Function Theorem, we see that
(\ref{betaequ}) has a solution for $\beta\in
(\beta_0-\delta,\beta_0+\delta)$ for some $\delta>0$. We denote it
by $\Sigma_{\beta}$.

It remains to show that $\Sigma_{\beta}$ is strictly stable. This
follows from (\ref{2nd2}). By (\ref{2nd2}), we see that the
operator $J_{\beta}$ is continuous in $\beta$. In particular, the
eigenvalues of $J_{\beta}$ are continuous in $\beta$. Since
$\Sigma_{\beta_0}$ is strictly stable, the smallest eigenvalue of
$J_{\beta_0}$ is positive. By continuity, there exists
$\delta_1>0$ (which may be smaller than $\delta$ above), so that
the smallest eigenvalue of $J_{\beta}$ is positive for $\beta\in
(\beta_0-\delta_1,\beta_0+\delta_1)$. That means, $\Sigma_{\beta}$
is a strictly stable $\beta$-symplecitc critical surface in $M$
for $\beta\in (\beta_0-\delta_1,\beta_0+\delta_1)$, namely,
$(\beta_0-\delta_1,\beta_0+\delta_1)\subset S$. This proves the
openness of $S$. \hfill Q.E.D.

\vspace{.1in}

\section{Rotationally Symmetric $\beta$-Symplectic Critical Surfaces in ${\mathbb C}^2$}

In this section, we study the rotationally symmetric
$\beta$-symplectic critical surfaces in ${\mathbb C}^2$ of the following form:
\begin{equation}\label{e5.1}
    F(r,\theta)=(r\cos\theta,r\sin\theta,f(r), g(r)).
\end{equation}
Then we have
\begin{equation*}
e_1=F_r=(\cos \theta, \sin \theta, f', g'), \ \ \ \
e_2=F_{\theta}=(-r\sin \theta, r\cos \theta, 0, 0),
\end{equation*}
\begin{equation*}
 e_3=(-f'\cos \theta, -f'\sin \theta, 1, 0), \ \ e_4=(-g'\cos\theta, -g'\sin\theta, 0, 1).
\end{equation*}
Further, we have
$$F_{rr}=(0, 0, f'', g''), \ F_{r\theta}=(-\sin \theta, \cos \theta, 0, 0), \ F_{\theta\theta}=(-r\cos \theta, -r\sin \theta, 0, 0).$$
Therefore, the induced metric on $\Sigma$ is given by
\begin{equation*}
(g_{ij})_{1\leq i,j, \leq 2}= \left(
  \begin{array}{cc}
    1+(f')^2+(g')^2 & 0 \\
        0       & r^2 \\
  \end{array}
\right),
\end{equation*}
\begin{equation*}
\det(g_{ij})=r^2(1+(f')^2+(g')^2),
\end{equation*}
\begin{equation*}
(g^{ij})_{1\leq i,j, \leq 2}= \left(
  \begin{array}{cc}
    \frac{1}{1+(f')^2+{(g')^2}} & 0 \\
        0       & \frac{1}{r^2} \\
  \end{array}
\right),
\end{equation*}
\begin{equation*}
(g_{\alpha\beta})_{3\leq \alpha,\beta \leq 4}= \left(
  \begin{array}{cc}
    1+(f')^2  &  f'g' \\
      f'g'   & 1+(g')^2 \\
  \end{array}
\right),
\end{equation*}
\begin{equation*}
\det(g_{\alpha\beta})=1+(f')^2+(g')^2.
\end{equation*}
\begin{equation*}
(g^{\alpha\beta})_{3\leq \alpha,\beta \leq 4}= \frac{1}{1+(f')^2+(g')^2}\left(
  \begin{array}{cc}
    1+(g')^2  &  -f'g' \\
      -f'g'   & 1+(f')^2 \\
  \end{array}
\right),
\end{equation*}
Denote by $A=1+(f')^2+(g')^2$, then it is easy to see that
\begin{equation}\label{e.cos}
    \cos\alpha=\frac{1}{\sqrt{1+(f')^2+(g')^2}}=\frac{1}{\sqrt{A}}>0,
\end{equation}
so that the surface $\Sigma$ is always symplectic. Now we compute ${\textbf H}$ and $(J(J\nabla\cos\alpha)^T)^\perp$.
\begin{eqnarray}\label{H}
{\textbf H} &=&g^{\alpha\beta}\langle g^{ij}\frac{\partial^2 F}{\partial x_i\partial x_j}, v_\alpha\rangle e_\beta\nonumber\\
&=&[g^{33}\langle g^{11}\frac{\partial^2 F}{\partial r^2}+g^{22}\frac{\partial^2 F}{\partial\theta^2},e_3\rangle
+g^{34}\langle g^{11}\frac{\partial^2 F}{\partial r^2}+g^{22}\frac{\partial^2 F}{\partial\theta^2},e_4\rangle]e_3\nonumber\\
&&+[g^{34}\langle g^{11}\frac{\partial^2 F}{\partial r^2}+g^{22}\frac{\partial^2 F}{\partial\theta^2},e_3\rangle
+g^{44}\langle g^{11}\frac{\partial^2 F}{\partial r^2}+g^{22}\frac{\partial^2 F}{\partial\theta^2},e_4\rangle]e_4\nonumber\\
&=&[\frac{1+(g')^2}{A}(\frac{f''}{A}+\frac{f'}{r})-\frac{f'g'}{A}(\frac{g''}{A}+\frac{g'}{r})]e_3\nonumber\\
&&+[\frac{1+(f')^2}{A}(\frac{g''}{A}+\frac{g'}{r})-\frac{f'g'}{A}(\frac{f''}{A}+\frac{f'}{r})]e_4\nonumber\\
&=&\frac{1}{rA^2}[r(1+(g')^2)f''-rf'g'g''+Af']e_3\nonumber\\&&+\frac{1}{rA^2}[r(1+(f')^2)g''-rf'g'f''+Ag']e_4.
\end{eqnarray}
It is clear that
\begin{eqnarray*}
\nabla\cos\alpha &=&(g^{11}\frac{\partial\cos\alpha}{\partial r}+g^{12}\frac{\partial\cos\alpha}{\partial\theta})e_1
+(g^{21}\frac{\partial\cos\alpha}{\partial r}+g^{22}\frac{\partial\cos\alpha}{\partial\theta})e_2\\
&=&g^{11}\frac{\partial\cos\alpha}{\partial r}e_1=-\frac{1}{A^{5/2}}(f'f''+g'g'')e_1.
\end{eqnarray*}
Moreover,
\begin{eqnarray*}
(Je_1)^T=g^{ij}\langle Je_1, e_i\rangle e_j=g^{22}\langle Je_1, e_2\rangle e_2=\frac{1}{r} e_2,
\end{eqnarray*}
\begin{eqnarray*}
(J(Je_1)^T)^\perp &=&\frac{1}{r}(Je_2)^\perp=\frac{1}{r}[g^{\alpha\beta}\langle Je_2, e_\alpha\rangle ]e_\beta\\
&=&\frac{1}{r}[\frac{1+(g')^2}{A}rf'-\frac{f'g'}{A}rg']e_3
  +\frac{1}{r}[\frac{1+(f')^2}{A}rg'-\frac{f'g'}{A}rf']e_4\\
&=&\frac{f'}{A}e_3+\frac{g'}{A}e_4.
\end{eqnarray*}
Thus,
\begin{eqnarray}\label{j1}
(J(J\nabla\cos\alpha)^T)^\perp &=&-\frac{(f'f''+g'g'')f'}{A^{7/2}}e_3-\frac{(f'f''+g'g'')g'}{A^{7/2}}e_4.
\end{eqnarray}
Putting (\ref{H}) and (\ref{j1}) into (\ref{betaequ})
\begin{eqnarray*}
 \left\{
  \begin{array}{cc}
    \frac{1}{rA^{7/2}}[r(1+(g')^2)f''-rf'g'g''+Af']  & = -\beta\frac{f'f''+g'g''}{A^{7/2}}f' \\
      \frac{1}{rA^{7/2}}[r(1+(f')^2)g''-rf'g'f''+Ag']  &  =-\beta\frac{f'f''+g'g''}{A^{7/2}}g'
  \end{array}
\right..
\end{eqnarray*}
That is equivalent to
\begin{eqnarray*}
 \left\{
  \begin{array}{cc}
    r(1+(g')^2+\beta (f')^2)f''+r(\beta-1)f'g'g''+(1+(f')^2+(g')^2)f'  & = 0 \\
     r(1+(f')^2+\beta (g')^2)g''+r(\beta-1)f'g'f''+(1+(f')^2+(g')^2)g'  & = 0
  \end{array}
\right..
\end{eqnarray*}
This implies a nice equation
\begin{eqnarray*}
 \left\{
  \begin{array}{cc}
   (rf'(1+(f')^2+(g')^2)^{\frac{\beta-1}{2}})' & = 0 \\
     (rg'(1+(f')^2+(g')^2)^{\frac{\beta-1}{2}})' & = 0
  \end{array}
\right..
\end{eqnarray*}
It is equivalent that
\begin{eqnarray}\label{rot1}
 \left\{
  \begin{array}{cc}
   rf'(1+(f')^2+(g')^2)^{\frac{\beta-1}{2}} & = C_1 \\
     rg'(1+(f')^2+(g')^2)^{\frac{\beta-1}{2}} & = C_2
  \end{array}
\right.,
\end{eqnarray} for any constant $C_1, C_2$.

There is a trivial case for the system (\ref{rot1}). Actually, if $C_1=C_2=0$, then $f, g$ must be constants and a $\beta$-symplectic critical surface must be a plane.

Next, we will divide four cases to study the solutions of (\ref{rot1}). The easiest one is the case that $\beta=1$, which follows directly from (\ref{rot1}):

\begin{theorem}\label{beta=1}
Let $\Sigma$ be a $1$-symplectic critical surface (i.e., $\beta=1$) in ${\mathbb C}^2$ of the form (\ref{e5.1}), then there exists constants $C_1$, $C_2$, $C_3$ and $C_4$, such that
\begin{equation}
f=C_1\ln r + C_3,  \ \ \ \ g=C_2\ln r+C_4.
\end{equation}
\end{theorem}

When $\beta=0$, i.e., $\Sigma$ is a minimal surface in ${\mathbb C}^2$, we have

\begin{theorem}\label{beta=0}
Let $\Sigma$ be a minimal surface (i.e., $\beta=0$) in ${\mathbb C}^2$ of the form (\ref{e5.1}). If $g\equiv constant$, we suppose $C_1=1$ in (\ref{rot1}). Then
\begin{equation*}
f(r)=\pm\cosh^{-1} r+C_3,
\end{equation*}
and $\Sigma$ is the catenoid.

If both $f$ and $g$ are not constants, we suppose $C_1=C_2=1$ in (\ref{rot1}). Then there exist constants $C_3$ and $C_4$, such that
\begin{equation}
f=\pm\cosh^{-1}\frac{r}{\sqrt 2}+ C_3,  \ \ \ \ g=\pm\cosh^{-1}\frac{r}{\sqrt 2}+C_4.
\end{equation}
\end{theorem}

\textbf{Proof:} If $g\equiv constant$, then the first equation of (\ref{rot1}) for $\beta=0$ becomes
\begin{equation*}
rf'(1+(f')^2)^{-\frac{1}{2}} = 1.
\end{equation*}
Rearranging it, we obtain that
\begin{equation*}
(f')^2=\frac{1}{r^2-1}.
\end{equation*}
Integrating gives the desired conclusion. The second case follows in the same way since we have $f'=g'$ in this case.
\hfill Q.E.D.

\vspace{.1in}

We will analyze the asymptotic behavior of the solutions.
Note that, if one of $C_1$, $C_2$ is zero, say, $C_2=0$, then from the second equality of (\ref{rot1}), we see that $g'\equiv0$. In this case, the first equation in (\ref{rot1}) becomes a nonlinear equation for $f'$. We can analyze it in the same way as we will do in the following theorems in which we deal with the general case.

Now we assume $\beta\neq1,0$ and $C_1\neq 0, C_2\neq 0$ in the following. The constants $C_1, C_2$ are determined by the initial data of $f'$ and $g'$. Without loss of generality, we assume $C_1= C_2=1$.

\begin{theorem}\label{thm.rot}
For any $\beta>0$ and $\varepsilon>0$, the equations
\begin{eqnarray*}
 \left\{
  \begin{array}{cc}
   rf'(1+(f')^2+(g')^2)^{\frac{\beta-1}{2}} & = 1, \\
     rg'(1+(f')^2+(g')^2)^{\frac{\beta-1}{2}} & = 1,\\
     f(\varepsilon) &=f_0,\\
     g(\varepsilon) &=g_0,
  \end{array}
\right.
\end{eqnarray*}
 have a unique $C^\infty$-solution on $[\varepsilon, +\infty)$. Moreover, $f'=g'$. As $r\to\infty$, we have the asymptotic expansion,
 $$f'=\frac{1}{r}-\frac{\beta-1}{r^3}+o(r^{-3}).$$
 As $r\to 0$, we have the asymptotic expansion,
 $$f'=2^{\frac{1-\beta}{2\beta}}r^{-\frac{1}{\beta}}-\frac{\beta-1}{\beta}2^{-\frac{3\beta+1}{2\beta}}r^{\frac{1}{\beta}}+o(r^{\frac{1}{\beta}}).$$
\end{theorem}

\textbf{Proof:}
 First we show that a solution exists on $[\varepsilon, +\infty)$.  Comparing the above equations, we get that $\frac{f'}{g'}=1$, i.e, $g'=f'$. Putting $f'=g'$ into the second equation we get that $rf'(1+2(f')^2)^{\frac{\beta-1}{2}}=1$. In the interval $[\varepsilon, +\infty)$, $f'$ is bounded. Thus there exists a solution in $[\varepsilon, +\infty)$.

Now we turn to the uniqueness. Assume there is a another solution $(\tilde{f}, \tilde{g})$ which satisfies the above equations. By the same argument, we have ${\tilde f}'={\tilde g}'$ and
 \begin{equation}\label{uniqueness}
 f'(1+2(f')^2)^{\frac{\beta-1}{2}} -{\tilde f}'(1+2({\tilde f}')^2)^{\frac{\beta-1}{2}}=0.
\end{equation}
 If $f'>{\tilde f}'$ at some point $r_0$, then $(1+2(f')^2)^{\frac{\beta-1}{2}}\geq(1+2({\tilde f}')^2)^{\frac{\beta-1}{2}}>0$ at $r_0$ when $\beta\geq1$. That is impossible. Therefore $f'\equiv {\tilde f}'$. Then $f={\tilde f}$ if they satisfy the same initial condition. Similarly,  $g={\tilde g}$. This proves the uniqueness for $\beta\geq1$.

On the other hand, if $0<\beta<1$, then we can rewrite (\ref{uniqueness}) as
\begin{equation*}
 f'(1+2({\tilde f}')^2)^{\frac{1-\beta}{2}} -{\tilde f}'(1+2(f')^2)^{\frac{1-\beta}{2}}=0.
\end{equation*}
Since $f'>0$, ${\tilde f}'>0$ by the equation, we can rearrange it as
\begin{equation*}
(f')^{\frac{2}{1-\beta}}-({\tilde f}')^{\frac{2}{1-\beta}}+2(f')^2({\tilde f}')^2 [(f')^{\frac{2\beta}{1-\beta}}-({\tilde f}')^{\frac{2\beta}{1-\beta}}]=0.
\end{equation*}
If $f'>{\tilde f}'$ at some point $r_0$, then $(f')^{\frac{2}{1-\beta}}>({\tilde f}')^{\frac{2}{1-\beta}}$ and $(f')^{\frac{2\beta}{1-\beta}}>({\tilde f}')^{\frac{2\beta}{1-\beta}}$ at $r_0$ when $0<\beta<1$. That is impossible. Therefore $f'\equiv {\tilde f}'$. Then $f={\tilde f}$ if they satisfy the same initial condition. Similarly,  $g={\tilde g}$. This proves the uniqueness for $0<\beta<1$.

Next we study the behaviour of $f'$ which satisfies
\begin{equation}\label{rot2} rf'(1+2(f')^2)^{\frac{\beta-1}{2}}=1,\end{equation}
as $r\to \infty$ and $r\to 0$.

As $r\to\infty$, $f'$ must tend to $0$. By the equation (\ref{rot2}), we see that $\lim_{r\to\infty} rf'=1$. We write $f'=\frac{1}{r}+\psi(r)$, then we obtain that $\lim_{r\to \infty} r\psi(r)=0$. Putting it into (\ref{rot2}),
$$(1+r\psi)^{\frac{2}{\beta-1}}(1+\frac{2}{r^2}+\frac{4\psi}{r}+2\psi^2)=1.$$
That implies that
$$\frac{2}{r^2\psi}+\frac{2}{\beta-1}r+O(\frac{1}{r})=0.$$ Therefore,
$$r^3\psi=-(\beta-1)+o(1),$$ which implies the desired expansion.

As $r\to 0$,  $f'$ must tend to infinity. By the equation (\ref{rot2}), we can easily see that $\lim_{r\to 0} 2^{\frac{\beta-1}{2}}r(f')^{\beta}=1$, which implies that $\lim_{r\to 0} r^{\frac{1}{\beta}}f'=2^{\frac{1-\beta}{2\beta}}$. Now we can write
$r^{\frac{1}{\beta}}f'=2^{\frac{1-\beta}{2\beta}}+w(r)$, where $\lim_{r\to 0}w(r)=0$. Putting it into (\ref{rot2}), we obtain that
$$(1+2^{\frac{\beta-1}{2\beta}}w(r))^{\frac{2}{\beta-1}}(2^{\frac{1}{\beta}}r^{\frac{2}{\beta}}+2^{\frac{\beta-1}{\beta}}(2^{\frac{1-\beta}{2\beta}}+w(r))^2)=1.$$ Expanding this equation we obtain that,
$$2^{-\frac{1}{\beta}}r^{\frac{2}{\beta}}+\frac{\beta}{\beta-1}w(r)+o(w(r))=0,$$ which implies that
$$\lim_{r\to 0}\frac{w(r)}{r^{\frac{2}{\beta}}}=-\frac{\beta-1}{\beta}2^{-\frac{3\beta+1}{2\beta}}.$$ Therefore, as $r\to 0$ we have,
$$f'=2^{\frac{1-\beta}{2\beta}}r^{-\frac{1}{\beta}}-\frac{\beta-1}{\beta}2^{-\frac{3\beta+1}{2\beta}}r^{\frac{1}{\beta}}+o(r^{\frac{1}{\beta}}).$$ \hfill Q.E.D.

\vspace{.1in}

The above theorem describes the asymptotic behavior of the rotationally symmetric $\beta$-symplectic critical surface as $r\to0$ and $r\to\infty$ for each fixed $\beta$. Next we will examine the behavior on a fixed compact set when $\beta$ goes to infinity and $\beta$ tends to 0.

\vspace{.1in}

\begin{proposition}
$f'\to 0$ as $\beta\to\infty$ in $[\varepsilon, +\infty)$ for any $\varepsilon>0$.
\end{proposition}

\textbf{Proof:}
Using the fundamental inequality $(1+x)^p\geq 1+px, $ for $x>0$ and $p>0$, we have that $(1+2(f')^2)^{\frac{\beta-1}{2}}\geq1+(\beta-1)(f')^2$ for $\beta>1$. By (\ref{rot2}), we have
$f'\leq \sqrt[3]{\frac{1}{(\beta-1)r}}$. This clearly implies the theorem.
\hfill Q.E.D.

\begin{corollary}
Let $\Sigma_{\beta}$ be a family of smooth complete $\beta$-symplectic critical surface which are rotationally symmetric. Suppose $\Sigma_{\beta}\cap K\neq\emptyset$ for some compact set $K\subset{\mathbb C}^2$ for all $\beta$. Then there exists a subsequence $\Sigma_{\beta_i}$, such that $\Sigma_{\beta_i}$ converges to a plane locally on ${\mathbb R}^4\backslash \{0\}$.
\end{corollary}

\textbf{Proof:}
For any $[a,b]\subset (0,\infty)$, we see from Theorem \ref{thm.rot} that $f'$ converges to $0$ uniformly on $[a,b]$ as $\beta\to\infty$. Therefore, by our assumption, there is a subsequence $\Sigma_{\beta_i}$, which converges to a piece of plane uniformly on $[a,b]$ as $\beta_i$ tends to infinity. Then the corollary follows from diagonal argument.
\hfill Q.E.D.

\vspace{.1in}

Now we study the limiting surface as $\beta\to0$.

\vspace{.1in}

\begin{proposition}\label{prop-beta-0}
For any $B>A>\sqrt2$, the solution $f_{\beta}$ of (\ref{rot2})
with $\beta>0$ converges uniformly on $[A,B]$ to the catenoid when
$\beta$ goes to zero. And for each fixed $r_0\in (0,\sqrt2]$, we
have $\lim_{\beta\to 0}f'_\beta(r_0)=\infty$.
\end{proposition}

\textbf{Proof:} We can rewrite (\ref{rot2}) as
\begin{equation}\label{rot3}
\frac{(f_{\beta}')^2}{(1+2(f_{\beta}')^2)^{1-\beta}}=\frac{1}{r^2}.
\end{equation}
On the interval $[A,B]$, we have $\frac{1}{r^2}\leq\frac{1}{A^2}<\frac{1}{2}$. Using the fundamental inequality $(1+x)^p\leq 1+x^p, $ for $x>0$ and $0<p<1$, we know from (\ref{rot3}) and Young's inequality that
\begin{equation*}
A^2(f_{\beta}')^2\leq 1+2^{1-\beta}(f_{\beta}')^{2(1-\beta)}\leq1+2(f_{\beta}')^{2(1-\beta)}\leq 2(1-\beta)(f_{\beta}')^2+1+2\beta.
\end{equation*}
In particular, for $\beta\in [0,1)$, we have
$$(f'_{\beta})^2\leq\frac{3}{A^2-2}.$$
Therefore, we can letting $\beta$ goes to zero in (\ref{rot3}) on $[A,B]$ and the first conclusion follows from the continuous dependence of ODE on the parameters.

For the second conclusion, we note that for each fixed $r_0\in (0,\sqrt2)$, we have from (\ref{rot3}) that
\begin{equation}\label{rot4}
r_0^2(f_{\beta}')^2(r_0)\geq(1+2(f_{\beta}'(r_0))^2)^{1-\beta}.
\end{equation}
Suppose $\liminf_{\beta\to0}f_{\beta}'(r_0)=D<\infty$, then there exists a sequence $\beta_i\to0$, such that $\lim_{\beta_i\to0}f_{\beta_i}'(r_0)=D$. Replace $\beta$ by $\beta_i$ in (\ref{rot4}) and letting $i\to\infty$ yields
\begin{equation*}
0\geq 1+(2-r_0^2)D^2\geq 1 \ \ \ \forall \ fixed \ r_0\in (0,\sqrt2],
\end{equation*}
which gives the desired contradiction.
\hfill Q.E.D.

\vspace{.1in}

In all the previous analysis, we choose $C_1=C_2=1$ in
(\ref{rot1}). Of course, we can also choose $C_1=C_2=-1$.
Actually, a smooth complete catenoid consists of two pieces of
graphs with $C_1=C_2=1$ (denoted by $\Sigma_0^+$) and with
$C_1=C_2=-1$ (denoted by $\Sigma_0^-$). By Theorem \ref{thm.rot},
we see that either $C_1=C_2=1$ or $C_1=C_2=-1$ produces a complete
$\beta$-symplectic critical surface $\Sigma_{\beta}^+$ (or
$\Sigma_{\beta}^-$ ).

Proposition \ref{prop-beta-0} tells us that, the two families of complete $\beta$-symplectic critical surfaces $\Sigma_{\beta}^+$ (or $\Sigma_{\beta}^-$) converges locally uniformly to the half of catenoid $\Sigma^+$ (or $\Sigma^-$, respectively) when $\beta\to0$. (This can also be seen from the Figure 1 in the following.)

Another interesting observation is the small perturbation
$\Sigma_{\beta}$ of the catenoid $\Sigma_0$ for $\beta>0$ consists
of two complete $\beta$-symplectic critical surfaces
$\Sigma_{\beta}^+$ and $\Sigma_{\beta}^-$. Note that the catenoid
$\Sigma_0$ is not a stable minimal surface, so one may get two
$\beta$-symplectic critical surfaces in a neighborhood of
$\Sigma_0$. It is different form the case that in Theorem
\ref{thm5.1}, where we show that there is only one stable
$\beta$-symplectic critical surface in a neighborhood of a stable
minimal surface.

  \begin{figure}[!htbp]
    \centering
    \includegraphics[width=.9\textwidth]{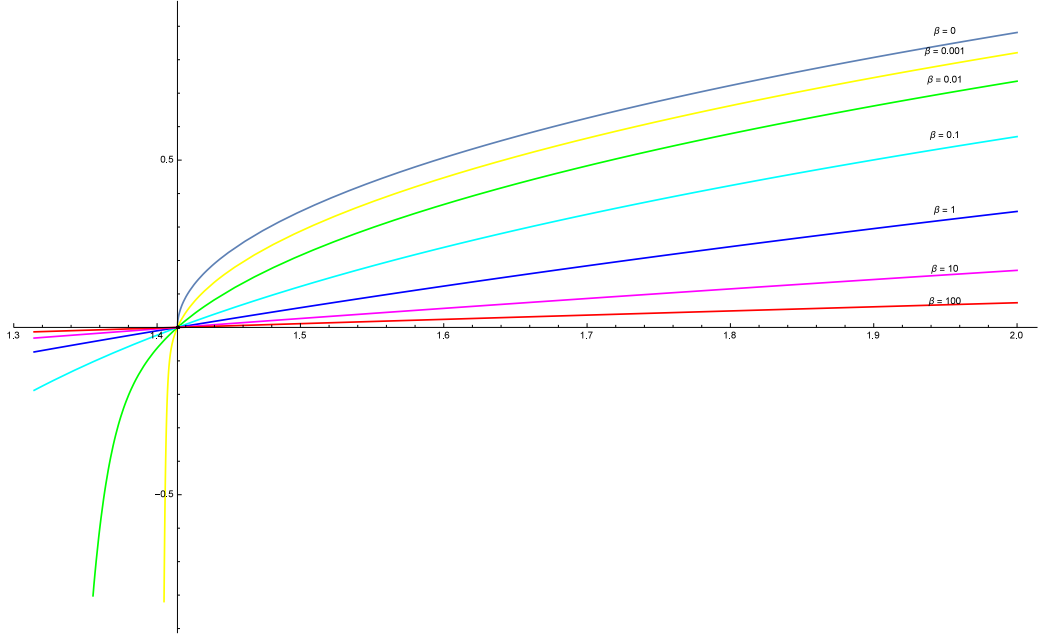}
    \caption{Comparison of different values of $\beta$}
  \end{figure}
\vspace{.1in}

\begin{corollary}
Let $\Sigma$ be a complete rotationally symmetric graphic $\beta$-symplectic critical surface in ${\mathbb C}^2$ with $\cos\alpha\geq \delta>0$ for some positive number $\delta$ and $\beta>0$, then $\Sigma$ must be a plane.
\end{corollary}

\textbf{Proof:} Theorem \ref{beta=1} and Theorem \ref{thm.rot} show that the domain of the defining functions $(f(r),g(r))$ must contain $(0,\infty)$. Furthermore, the asymptotic expansion show that if $\Sigma$ is not a plane, then either $\lim_{r\to 0}|f'(r)|=\infty$ or $\lim_{r\to 0}|g'(r)|=\infty$. By (\ref{e.cos}), we must have $\lim_{r\to 0}\cos\alpha=0$.
\hfill Q.E.D.

\vspace{.1in}

Finally, we prove a Louville theorem for $\beta$-symplectic critical surfaces in ${\mathbb C}^2$.

\begin{theorem}
Let $\Sigma^2$ be a complete $\beta$-symplecctic surface in ${\mathbb C}^2$ with $\cos\alpha\geq\delta>0$ and $Area(B^{\Sigma}(s))\leq Cs^2$ for $s>0$, then $\Sigma$ is a holomorphic curve with respect to some compatible complex structure in ${\mathbb C}^2$. Here, $B^{\Sigma}(s)$ is the intrinsic ball of $\Sigma$.
\end{theorem}

\textbf{Proof:}
From Corollary \ref{cor.angle}, we can easily see that:
\begin{eqnarray*}
\Delta \frac{1}{\cos\alpha}
&=& -\frac{\Delta\cos\alpha}{\cos^2\alpha}+\frac{2|\nabla\cos\alpha|^2}{\cos^3\alpha}\nonumber\\
&=&\frac{2\cos^4\alpha-2\beta\sin^4\alpha}{\cos^3\alpha
(\cos^2\alpha+\beta\sin^2\alpha)}|\nabla\alpha|^2
+\frac{2\sin^2\alpha}{\cos^3\alpha}|\nabla\alpha|^2\nonumber\\
&=&\frac{2}{\cos\alpha
(\cos^2\alpha+\beta\sin^2\alpha)}|\nabla\alpha|^2.
\end{eqnarray*}
Since $\cos\alpha\geq\delta>0$, we have $1\leq\frac{1}{\cos\alpha}\leq\frac{1}{\delta}$, which means that $\frac{1}{\cos\alpha}$ is a subharmonic function bounded from above on $\Sigma$. However, the quadratic area growth implies that $\Sigma$ is parabolic (\cite{CY}). This forces $\frac{1}{\cos\alpha}$ to be a constant on $\Sigma$, which implies that $\Sigma$ is holomorphic with respect to some compatible complex structure in ${\mathbb C}^2$.
\hfill Q.E.D.

\end{document}